\documentclass[12pt]{article}
\newcommand{\lb}[1]{\label{#1}}\newcommand{\reff}{\ref}\newcommand{\rf}[1]{(\reff{#1})}\newcommand{\cit}[1]{\cite{#1}}\newcommand{\bibi}[1]{\bibitem{#1}}%pour ne pas indiquer labels, citations et r\'ef\'erences
\newcommand{\com}[1]{}%_Pour ne pas indiquer les commentaires
\usepackage{amssymb}
\usepackage{amsmath}
\usepackage[latin1]{inputenc}
\usepackage{epsf}

\newtheorem{lm}{{\bf Lemma}}[section]
\newtheorem{theor}[lm]{{\bf Theorem}}
\newtheorem{deff}[lm]{{\bf Definition}}

\newtheorem{cj}[lm]{{\bf Conjecture}}

\newtheorem{cor}[lm]{{\bf Corollary}}
\newtheorem{exemple}[lm]{{\bf Example}}

\newtheorem{prop}[lm]{{\bf Proposition}}
\newcommand{\sep}{$\!${\rm.}\ }
\newcommand{\df}[2]{\begin{deff}{\lb{#1}}\sep{\rm #2}\end{deff}}

\newcommand{\lem}[2]{\begin{lm}{\lb{#1}}\sep{\sl #2}\end{lm}}
\newcommand{\theo}[2]{\begin{theor}{\lb{#1}}\sep{\sl 
#2}\end{theor}}
\newcommand{\conj}[2]{\begin{cj}{\lb{#1}}\sep{\sl #2}\end{cj}}
\newcommand{\coro}[2]{\begin{cor}{\lb{#1}}\sep{\sl #2}\end{cor}}
\newcommand{\propo}[2]{\begin{prop}{\lb{#1}}\sep{\sl #2}\end{prop}}
\renewcommand{\sec}[2]{\section{#2}{\lb{#1}}}

\newcommand{\eq}[2]{\begin{equation}#2\ \ \  \ \lb{#1}\end{equation}}

\newcommand{\pr}{\noindent {\sl Proof\;\sep }}
\newcommand{\ep}{\hfill \framebox[2mm]{\ } \medskip}
\newcommand{\tq}{\ ;\ }
\newcommand{\be}{\begin{enumerate}}
\newcommand{\ee}{\end{enumerate}}
\newcommand{\noi}{\noindent}
\newcommand{\med}{\medskip}

\newcommand{\ds}{\displaystyle}
\newcommand{\lp}{\left(}
\newcommand{\rp}{\right)}

\newcommand{\dist}{\,\mbox{\rm dist}}
\newcommand{\lng}{\,\mbox{\rm length}}

\renewcommand{\hat}{\widehat}

\newcommand{\ssi}{\Leftrightarrow}
\renewcommand{\a}{\alpha}

\renewcommand{\b}{\beta}
\newcommand{\cc}{{\bf C}}
\renewcommand{\d}{\delta}

\newcommand{\De}{\Delta}

\newcommand{\eps}{\varepsilon}

\newcommand{\g}{\gamma}

\renewcommand{\O}{{\mathcal O}}
\renewcommand{\o}{{\mathit o}}
\newcommand{\q}{\theta}

\renewcommand{\S}{{\mathbb S}}

\newcommand{\N}{\mathbb{N}}
\newcommand{\Z}{\mathbb{Z}}
\newcommand{\Q}{\mathbb{Q}}
\newcommand{\R}{\mathbb{R}}
\newcommand{\C}{\mathcal{C}}

\newcommand{\nt}{{\rm Int}\,}

\newcommand{\conv}{{\rm conv}}
\newcommand{\cotan}{{\rm cotan}}

\newcommand{\ed}{\end{document}}

\newcommand{\vecc}[1]{\stackrel{\rightarrow}{#1}}
\newcommand{\veccc}[1]{\stackrel{\longrightarrow}{#1}}
\title{Envelopes of $\alpha$-sections}
\author{N.\ Chevallier, A.\ Fruchard, and C.\ V\^ilcu}
\date{\empty}
\begin{document}
\maketitle
\abstract{
\noi Let $K$ be a planar convex body af area $|K|$, and take $0<\a<1$.
An $\a$-section of $K$ is a line cutting $K$ into two parts, one of which
has area $\a|K|$. 
This article presents a systematic study of the envelope of $\a$-sections
and its dependence on $\a$. Several open questions are asked, one of them
in relation to a problem of fair partitioning.
}

\

{\noi\bf Keywords:} Convex body, alpha-section, envelope, floating body, 
%alpha-core, 
fair partitioning.

\

{\noi\bf MSC Classification: 52A10, 52A38, 51M25, 51M04}

\
%
%%%%%%%%%%%%%%%%%%%%%%%%%%%%%%%%%%%%%%%%%%%%%%%%
%
\sec{1.}{Introduction}
In this paper, unless explicitly stated otherwise, $K$ denotes a convex 
body in the Euclidean plane $E$; {\sl i.e.,} a compact convex subset of 
$E$ with nonempty interior. Let $\partial K$ denote the boundary of $K$ 
and $|K|$ its area.
Given $\a\in\,]0,1[$, an {\sl$\a$-section of $K$} is an oriented
line $\De \subset E$ cutting $K$ in two parts, one to the right, 
denoted by $K^-$, of area $|K^-|=\a|K|$,
and the other to the left, $K^+$, of area $|K^+|=(1-\a)|K|$; 
here $K^\pm$ are compact sets, thus $K^+\cap K^-=\De\cap K$.

Denote by $K_\a$ the intersection of all $K^+$ and call it the
{\sl$\a$-core} of $K$; 
denote by $m_\a$ the envelope of all $\a$-sections of $K$.

The purpose of this article is to study $m_\a$ and its relation to the 
$\a$-core. We refer to Sections~\reff{5.}-\reff{6.} for formal statements,
and give in this introductory section an informal presentation.
Since $K_\a$ is empty for $\a>\frac12$,  we will often implicitly assume 
$\a\leq\frac12$ when dealing with $\a$-cores.
The situation depends essentially upon whether $K$ is centrally symmetric 
(we will say `symmetric' in the sequel, for short) or not.
Precisely, we prove the following statements.

If $K$ is symmetric then one has $m_\a=\partial K_\a$ for all
$\a\in\,\big]0,\tfrac12\big[$.
Moreover, we have the following equivalence:
the envelope $m_\a$ is of class $\C^1$ for all $\a\in\,\big]0,\tfrac12\big[$
if and only if $K$ is strictly convex.

If $K$ is non-symmetric then we cannot have $m_\a=\partial K_\a$ for all 
$\a\in\,\big]0,\tfrac12\big[$, because 
$m_\a$ exists for all $\a$, whereas $K_\a$ is empty for $\a$ close enough 
to $\frac12$. More precisely, there exists a critical value
$\a_B\in\big[0,\tfrac12\big[$ such that 
for all $\a\in\,]0,\a_B]$ we have $m_\a=\partial K_\a$, 
and for all $\a\in\,\big]\a_B,\tfrac12\big[$ 
we have $m_\a\supsetneq\partial K_\a$.

The case $\a_B=0$ can occur, {\sl e.g.}, if there exists a triangle
containing $K$ with an edge entirely contained in $\partial K$. 
We also prove that $m_\a$ is never of class $\C^1$ for 
$\a\in\,\big]\a_B,\tfrac12\big[$, and that 
$m_\a$ is of class $\C^1$ for all $\a\in\,]0,\a_B[$ if and only if 
$\partial K$ does not contain two parallel segments. 
As a by-product, we obtain the following characterization: A convex 
body $K$ is non-symmetric if and only if there exists a triangle
containing more than half of $K$ (in area), with one side entirely 
in $K$ and the two others disjoint from the interior of $K$.

Concerning the $\a$-core, we prove that there is another critical value, 
$\a_K\in\big[\frac49,\frac12\big]$, such that 
if $0<\a<\a_K$ then $K_\a$ is strictly convex with nonempty interior,
if $\a=\a_K$ then $K_\a$ is  reduced to one point,
and if $\a_K<\a<1$ then $K_\a$ is empty.
We emphasize that, when $K_\a$ is a point, this point is 
{\sl not} necessarily the mass center of $K$, see Section 8.9. 
The value $\a_K=\frac12$ occurs if and only if $K$ is symmetric and the
value $\a_K=\frac49$ occurs if and only if $K$ is a triangle.

A similar study, for secants between parallel supporting lines to $K$, 
whose distances to the corresponding lines make a ratio of $\a/(1-\a)$,
is the subject of {\sl (ir)reducibility theory} of convex bodies.
There too, the envelope of those secants is sometimes different from the 
intersection of the half-planes they are defining; and there exists a ratio,
called critical, for which the later object is reduced to a point. 
See for example \cit{h,kl,z1}.

Our paper is closely related to previous works about slicing convex bodies,
outer billiards (also called {\sl dual billiard}), floating bodies,
and fair (or equi-) partitioning.
It is also related to continuous families of curves in the sense of
Gr\"unbaum, see Section 8.8 and the references therein.
There is a vast literature on these subjects; we refer in the following to 
very few articles, and briefly present even fewer, that we find particularly 
relevant for our study.
Further references can be found in those papers.

Generalizing previous results on common tangents and common 
transversals to families of convex bodies \cit{bhj}, \cit{cgppsw}, 
J.~Kincses~\cit{k} showed that, 
for any {\sl well-separated} family of strictly convex bodies,
the space of $\a$-sections is diffeomorphic to $\mathbb{S}^{d-k}$. 

A {\sl billiard table} is a planar strictly convex body  $K$. 
Choose a starting point $x$ outside the table and one of the two tangents
through $x$ to $K$, say the right one, denoted by $D$; the image $T(x)$ of
$x$ by the {\sl billiard map} $T$ is the point symmetric of $x$ with respect
to the tangency point  $D\cap\partial K$. 
A {\sl caustic} of the billiard is an invariant curve (an invariant torus
in the terminology of the KAM theory).
The link between outer billiards and $\a$-sections is the following: 
If $K$ is the envelope of $\a$-sections of a convex set bounded by a curve $L$,
for some $\a$, then $L$ is a caustic for the outer billiard of table $K$,
{\sl cf. e.g.}~\cit{ft}.

Outer billiards have been considered by several authors, such as 
J.~Moser~\cit{mo}, V.~F.~La\-zutkin~\cit{l},
E.~Gutkin and A.~Katok~\cit{gk}, D.~Fuchs and S.~Tabachnikov~\cit{ft}, 
S.~Tabachnikov~\cit{t1}.
Therefore, it is natural that several authors considered envelopes of 
$\a$-sections in the framework of outer billiard,
see {\sl e.g.}, Lecture~11 of the book of D.~Fuchs and 
S.~Tabachnikov~\cit{ft}, and references therein.
Besides those studies, the envelope of $\a$-sections seems to be scarcely
studied.

With our notation, the set $K_{[\a]}$ bounded by $m_\a$ was called 
{\sl floating body} of $K$ and its study goes back to C.~Dupin, 
see \cit{du,w}.
On the other hand, our ``$\a$-core'' $K_\a$ was introduced by C.~Sch\"utt
and E.~Werner \cit{sw} and studied in a series of papers 
\cit{sw,sw2, st, w}.
under the name of {\sl convex floating body}. For convex bodies $K$ in 
$\mathbb{R}^{d}$ and for $\a$ small enough, 
they gave estimates for
${\rm vol}_n(K)-{\rm vol}_n(K_{[\a]})$ and for
${\rm vol}_n(K)-{\rm vol}_n(K_\a)$, in relation to the affine surface area
and to polygonal approximations.
M.~Meyer and S.~Reisner~\cit{mr} proved in arbitrary dimension that 
$K$ is symmetric if and only if $m_\a=\partial K_\a$ for any 
$\a\in\,\big]0,\tfrac12\big[\,$. They also prove that $m_\a$ is smooth if 
$K$ is strictly convex. 
A. Stancu \cit{st} considered convex bodies $K \subset \mathbb{R}^{d}$
with boundary of class ${\cal C}^{\geq 4}$, and proved that there exists
$\delta_K >0$ such that $K_{\delta}$ is homothetic to $K$, for some
$\delta < \delta_K$, if and only $K$ is an ellipsoid.

The terminology of ``(convex) floating body'' is very suggestive for the
floating theory in mathematical physics. For our study however,
considering its close connections to $\a$-sections and fair partitioning,
it seems more natural to use the term of ``$\a$-core''.

Some of our results, essentially Proposition~\reff{p1},
are known and already published.
In that case we mention references after the statement.
For the sake of self-containedness, however, we will provide complete proofs.

We end the paper with several miscellaneous results and open questions.
Our main conjecture is as follows: If $K\subseteq L$ are two convex bodies
and $\a\in\,\big]0,\tfrac12\big[$, then there exists an $\a$-section of $L$
which either does not cross the interior of $K$, or is a $\b$-section of $K$
for some $\b\leq\a$. The conjecture has been recently proven in the case of 
planar convex bodies, see~\cit{fm}.

There are two reasons which motivated us to undertake this systematic 
study of $\a$-sections and their envelopes.
First, we found only one reference which focuses especially on the envelope
of $\a$-sections, Lecture 11 of the nice book \cit{ft},
which is a simplified approach.
Although the results contained in the present article seem natural, 
and the proofs use only elementary tools and are most of the time simple, 
we hope that our work will be helpful in clarifying the things.

Our second motivation for studying $\a$-sections and their 
envelope was a problem of fair partitioning of a pizza.
%which we present in Section~\reff{8.}.
What we call {\sl a pizza} is a pair of planar convex bodies $K\subseteq L$,
where $L$ represents the {\sl dough} and $K$ the {\sl topping} of the pizza. 
The problem of fair partitioning of convex bodies in $n$ pieces 
is a widely studied topic, see {\sl e.g.}~\cit{bbs,bm,b,bz,ka,kha,s,so}.
%we give a very short account of the literature in Section~\reff{8.}.
Nevertheless, to our knowledge, our way of cutting has never been considered:
We use a succession of double operations: a cut by a {\sl full}
straight line, followed by a Euclidean move of one of the resulting pieces;
then we repeat the procedure.
The final partition is said to be {\sl fair} if each resulting slice has 
the same amount of $K$ and the same amount of $L$. The result of~\cit{fm}
is the following: 
Given an integer $n\geq2$, there exists a fair partition of {\sl any} pizza 
$(K,L)$ into $n$ parts if and only if $n$ is even.
%
%
%%%%%%%%%%%%%%%%%%%%%%%%%%%%%%%%%%%%%%%%%%%%%%%%
%
\sec{2.}{Notation, conventions, and preliminaries}
The notation $\S^1$ stands for the standard unit circle, $\S^1:=\R/(2\pi\Z)$,
endowed with its usual metric
%$d(\q,\q')=\min(|\q-\q'|,|\q-\q'-2\pi|,|\q-\q'+2\pi|)$.
$d(\q,\q')=\min\{|\tau-\tau'|\tq\tau\in\q,\;\tau'\in\q'\}$.
On $\S^1$ we use the notation:
%natural ``local order'': 
$\q\leq\q'$ if there exist 
$\tau\in\q,\;\tau'\in\q'$ such that $\tau\leq\tau'<\tau+\pi$;
it is not an order because it is not transitive.

Given $\q\in\S^1$, let $\vec u(\q)$ denote the unit vector of direction $\q$,
$\vec u(\q)=(\cos\q,\sin\q)$.
For convenience, we add arrows $\,\vec{\;}\,$ on vectors.
Unless explicitly specified otherwise, all derivatives will be with respect to
$\q$, hence {\sl e.g.}, $\vec u\,'(\q)=\frac{d\vec u}{d\q}(\q)=(-\sin\q,\cos\q)$ 
is the unit vector orthogonal to $\vec u(\q)$ such that the frame 
$\big(\vec u(\q),\vec u\,'(\q)\big)$ is counterclockwise.

Given an oriented straight line $\De$ in the plane, $\De^+$ denotes the
{\sl closed} half-plane on the left bounded by $\De$, 
and $\De^-$ is the closed half-plane on the right.
We identify oriented straight lines with points of the cylinder 
$\cc=\S^1\times\R$, associating each pair $(\q,t)\in\cc$ to the line
oriented by $\vec u(\q)$ and passing at the signed distance $t$ from
the origin.
In other words, the half-plane $\De^+$ is given by 
$\De^+=\{x\in\R^2\tq\langle x,\vec u\,'(\q)\rangle\geq t\}$.
We endow $\cc$ with the natural distance 
$d\big((\q,t),(\q',t')\big)=\big(d(\q,\q')^2+|t-t'|^2)\big)^{1/2}$.
\med

Given $\a\in\,]0,1[$, an {\sl $\a$-section of $K$} is an oriented line
$\De$ such that $|\De^-\cap K|=\a|K|$.

For all $\a\in\,]0,1[$ and all $\q\in[0,2\pi[$, there exists a unique
$\a$-section of $K$ of direction $\q$; it will be denoted by $\De(\a,\q)$.
This defines a continuous function $\De:\,]0,1[\,\times\S^1\to\cc$.
We obviously have the symmetry
\eq d{
\De^\pm(1-\a,\q)=\De^\mp(\a,\q+\pi).
}
For $0<\a\leq\frac12$, we call {\sl $\a$-core of $K$}, and denote by $K_\a$,
the intersection of all left half-planes bounded by $\a$-sections:
$$
K_\a=\bigcap_{\q\in\S^1}\De^+(\a,\q).
$$
It is a compact convex subset of the plane, possibly reduced to one point
or empty. 
\med

Let $\De=\De(\a,\q)$ be an $\a$-section of $K$ of direction 
$\vec u=\vec u(\q)$, and let $b,c$ denote the endpoints of the chord
$\De\cap K$, with $\vec{bc}$ having the orientation of $\vec u$
({\sl i.e.}, the scalar product $\langle\vec{bc},\vec u\rangle$ is positive), 
see Figure~\reff{f1}.
\begin{figure}[htb]
\vspace{-7mm}
\begin{center}
\raisebox{0cm}{\epsfysize5cm\epsfbox{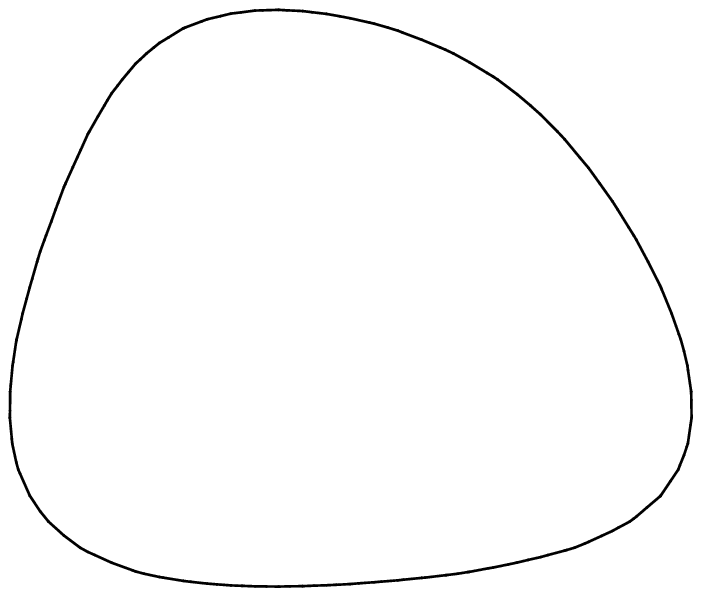}}

\vspace{-4cm}
\begin{picture}(300,100)
\put(60,0){\line(5,1){220}}\put(280,44){\vector(4,1){0}}
\put(40,60){\line(5,1){20}}\put(60,64){\vector(4,1){0}}
\put(50,50){$\vec u$}
\put(40,60){\line(-1,5){4}}\put(36,80){\vector(-1,4){0}}
\put(21,70){$\vec u\,'$}
\put(103,8.6){\line(5,-2){45}}\put(148,-9.4){\vector(3,-1){0}}
\put(103,8.6){\circle*{3}}\put(151,-13){$\vec b\,'$}
\put(100,-3){$b$}\put(125,4.5){$\b$}
\put(158.65,19.8){\circle*{3}}
\put(155,11){$m$}
\put(214.3,31){\vector(1,4){10}}\put(214.3,31){\circle*{3}}
\put(227,65){$\vec c\,'$}
\put(215,22){$c$}\put(222,40){$\g$}
\put(130,60){$K$}\put(180,100){$\partial K$}
\put(240,25){$\Delta$}%(\alpha,\theta)$}
\put(255,5){$\Delta^-$}\put(250,55){$\Delta^+$}
\end{picture}
\end{center}
\caption{{\small Some notation}}
\vspace{-5pt}
\lb{f1}
\end{figure}

Let $m=\tfrac12(b+c)$ denote the midpoint of the chord $bc$.
Let $h$ denote the half-length of $bc$; hence we have $b=m-h\vec u$
and $c=m+h\vec u$.
The functions $b,c,m$ and $h$ are continuous with respect to $\a$ and $\q$.
As we will see in the proof of Proposition~\reff{p1}, they are also left- and
right-differentiable with respect to $\q$ at each $\q_0\in\S^1$, {\sl e.g.},
the following limit exists
(with the convention  $\q\to\q_0^-$ for $\q\to\q_0,\q<\q_0$)
$$
\vec b_l'(\q_0)=\lim_{\q\to\q_0^-}
\tfrac1{\q-\q_0}\big(b(\a,\q)-b(\a,\q_0)\big),
$$
and similarly for $\vec b_r', \vec c_l\,', \vec c_r\,'$.
They are also left- and 
right-differentiable with respect to $\a$, but we will not use this fact.

We say that $b$ is {\sl a regular point of $\partial K$}
({\sl regular} for short) if there is a unique supporting line to $K$ at $b$
({\it i.e.}, $\vec b_l'=\vec b_r'$); otherwise we call $b$
{\sl a corner point of $\partial K$} (a {\sl corner} for short).

If $b$ is regular then $\b=(\hat{\vec b\,',\vec u})\in\,]0,\pi[$ denotes
the angle between the tangent to $\partial K$ at $b$ and $\De$.
If $b$ is a corner then $\b_l=(\hat{\vec b_l',\vec u})$, resp.
$\b_r=(\hat{\vec b_r',\vec u})$, is the angle between the left-tangent,
resp. right-tangent, to $\partial K$ at $b$ and $\De$.
Similarly, let $\g=(\hat{\vec u,\vec c\,'})\in\,]0,\pi[$ (if $c$ is regular),
resp. $\g_l,\,\g_r$ (if $c$ is a corner),
be the angle between $\De$ and the tangent, resp. left-tangent,
right-tangent, to $\partial K$ at $c$, see Figure~\reff{f1}.

For a fixed $\a$, the values of $\q$ such that $b(\a,\q)$ or $c(\a,\q)$
(or both) is a corner will be called {\sl singular}; those for which
both $b$ and $c$ are regular will be called {\sl regular}.

Observe that we always have $\b_r\leq\b_l$ and $\g_l\leq\g_r$,
with equality if and only if $b$, resp. $c$, is regular. 
Also observe the following fact.
\eq b{
b\mbox{ and }c\mbox{ admit parallel supporting lines if and only if }
\b_r+\g_l\leq\pi\leq\b_l+\g_r.
}
Finally, observe that the angle between $\vec b'$, (resp. $\vec b_l'$,
$\vec b_r'$) and the axis of abscissae 
is equal to $\q-\b(\a,\q)$ (resp. $\q-\b_l(\a,\q)$, $\q-\b_r(\a,\q)$).
Since these angles are increasing and intertwining functions of $\q$,
we have the following statement.
$$
\mbox{If }\q<\q'\mbox{ then }\q-\b_l(\a,\q)\leq\q-\b_r(\a,\q)
\leq\q'-\b_l(\a,\q')\leq\q'-\b_r(\a,\q').
$$
It follows that
\eq t{
\lim_{\q\to\q_0^+}\b_l(\a,\q)=\lim_{\q\to\q_0^+}\b_r(\a,\q)=
\b_r(\a,\q_0)\leq\b_l(\a,\q_0)
=\lim_{\q\to\q_0^-}\b_l(\a,\q)=\lim_{\q\to\q_0^-}\b_r(\a,\q),
}
and similarly
\eq{3b}{
\lim_{\q\to\q_0^-}\g_l(\a,\q)=\lim_{\q\to\q_0^-}\g_r(\a,\q)
=\g_l(\a,\q_0)\leq\g_r(\a,\q_0)
=\lim_{\q\to\q_0^+}\g_l(\a,\q)=\lim_{\q\to\q_0^+}\g_r(\a,\q).
}
In a same way we have
\eq g{
\mbox{If }\a<\a'\mbox{ then }\b_r(\a,\q)\leq\b_l(\a,\q)
\leq\b_r(\a',\q)\leq\b_l(\a',\q)
$$
$$
\mbox{ and }
\g_l(\a,\q)\leq\g_r(\a,\q)
\leq\g_l(\a',\q)\leq\g_r(\a',\q).
}
As a consequence,
$$
\lim_{\a\to\a_0^-}\b_l(\a,\q)=\lim_{\a\to\a_0^-}\b_r(\a,\q)=
\b_r(\a_0,\q)\leq\b_l(\a_0,\q)
=\lim_{\a\to\a_0^+}\b_l(\a,\q)=\lim_{\a\to\a_0^+}\b_r(\a,\q),
$$
$$
\lim_{\a\to\a_0^-}\g_l(\a,\q)=\lim_{\a\to\a_0^-}\g_r(\a,\q)
=\g_l(\a_0,\q)\leq\g_r(\a_0,\q)
=\lim_{\a\to\a_0^+}\g_l(\a,\q)=\lim_{\a\to\a_0^+}\g_r(\a,\q).
$$
All these elements $b,c,m,h,\b,\g$ can be considered as functions of both 
$\a$ and $\q$.
Nevertheless, as already said, all derivatives are with respect to $\q$.

Let $v=v(\a,\q)$ be the scalar product $v=\langle\vec m',\vec u\rangle\in\R$.
If $v$ has a discontinuity, then $v_l$ and $v_r$ denote its corresponding
left and right limits.
As we will see, $\vec m'$ is always collinear to $\vec u$, hence $v$ is
the ``signed norm'' of $\vec m'$.
We will also see that $v$ has discontinuities only if $b$ or $c$ (or both)
is a corner.
Since we chose $\q$ as parameter, we will also see that $v$ is the signed
radius of curvature of the curve $m$,
but we prefer to refer to it as the {\sl velocity} of the current point
$m$ of the envelope.
The symmetry~\rf d gives $m(\a,\q+\pi)=m(1-\a,\q)$, 
hence $\vec m'(\a,\q+\pi)=\vec m'(1-\a,\q)$.
Since $\vec u(\q+\pi)=-\vec u(\q)$, we obtain
\eq w{
v_l(\a,\q+\pi)=-v_l(1-\a,\q)\;\mbox{ and }\;v_r(\a,\q+\pi)=-v_r(1-\a,\q).
}
Our last notation is $V$ for the segment of endpoints $v_l$ and $v_r$;
{\sl i.e.}, $V=[v_l,v_r]$ if $v_l\leq v_r$, $V=[v_r,v_l]$ otherwise. 
Formulae~\rf w yield
\eq W{
V(\a,\q+\pi)=-V(1-\a,\q).
}
%
%%%%%%%%%%%%%%%%%%%%%%%%%%%%%%%%%%%%%%%%%%%%%%%%%%%%%%%%%%%%
%
\sec{2b.}{A ``digressing tour''}
Before carrying on $\a$-sections, we would like to digress for a moment
in a more general framework.
The notion and results of this section are elementary and probably 
already known, but we did not find any reference in the literature.
They can be considered as exercises in a graduate course on planar curves.

Recall that a {\sl ruled function} (or {\sl regulated function})
$R:\S^1\to\R$ is the uniform limit of piecewise constant functions. 
It is equivalent to saying that $R$ admits a left- and a right-limit,
denoted below by $R_l$ and $R_r$, at
each point of $\S^1$, see {\sl e.g.},~\cit{di}.
Recall also our notation $\vec u(\q)=(\cos\q,\sin\q)$.
\df{d1}{
{\noi\rm(a)} \ 
A {\sl tour} is a planar curve parametrized by its tangent. 
More precisely, we call $m:\S^1\to\R^2$ a {\sl tour} if $m$ is
continuous, has a left- and a right-derivative at each point of $\S^1$, 
and if there exists a ruled function $R:\S^1\to\R$ such that 
$\vec m'_l(\q)=R_l(\q)\vec u(\q)$ and $\vec m'_r(\q)=R_r(\q)\vec u(\q)$.
\med

{\noi\rm(b)} \ 
The {\sl core} $K=K(m)$ of a tour $m$ is the intersection of all left
half-planes delimited by all tangents to $m$ oriented by $\vec u$,
{\sl i.e.}, $K=\bigcap_{\q\in\S^1}D^+(\q)$, where 
$D(\q)=m(\q)+\R\vec u(\q)$. 
}
Tours are not necessarily simple curves. 
The case $m(\q)=m(\q+\pi)$ gives rise to a {\sl double half-tour}. 
This is the case {\sl e.g.} for the envelope of half-sections of a planar
convex body.
\begin{figure}[htb]
\begin{center}
\raisebox{0cm}{\epsfysize5cm\epsfbox{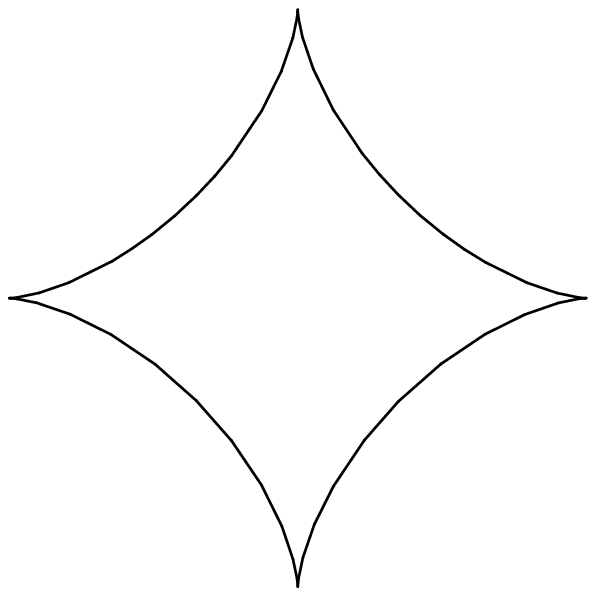}}
\hspace*{5cm}
\epsfysize5cm\epsfbox{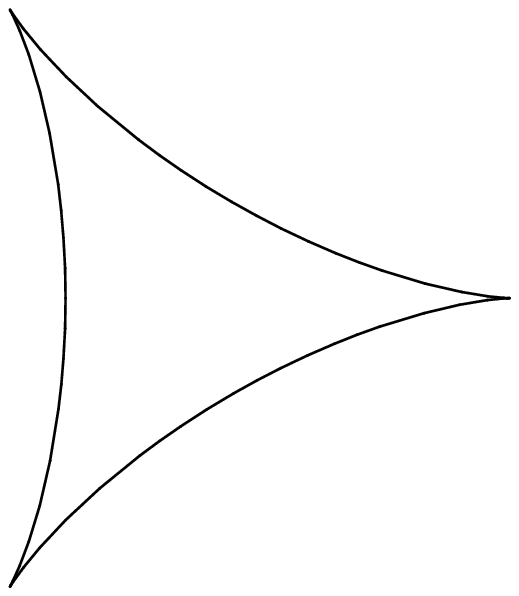}
\hspace*{17mm}
\end{center}

\vspace{-7mm}
\caption{{\small Some (double half-)tours: Left the astroid, right the deltoid.
The astroid cannot be an envelope of an $\a$-section of a planar convex 
body; the deltoid is probably such a half-section, although we ignore
how to prove it.}}
\vspace{-5pt}
\lb{f8}
\end{figure}
\propo{p1.1}{
Let $m$ be a tour, with associated ruled function $R$, and let $m^*$
denote its image in the plane: $m^*=m(\S^1)$.
\med

{\noi\rm(a)} \ 
For each $\q\in\S^1$, $R_l(\q)$, resp. $R_r(\q)$, is the signed
left-, resp. right-radius of curvature of $m^*$ at $m(\q)$.
\med

{\noi\rm(b)} \ 
If $R_l$ and $R_r$ do not vanish on $\S^1$, then $m^*$ is a
$\C^1$ submanifold of the plane.
\med

{\noi\rm(c)} \ 
Conversely, if there exist $\q_1\leq\q_2\leq\q_1+\pi\in\S^1$ such that
the product $R_r(\q_1)R_l(\q_2)$ is negative, then there exists
$\q_3\in[\q_1,\q_2]$ such that $m^*$ is not of class $\C^1$ at $m(\q_3)$.
\med

{\noi\rm(d)} \ 
The core $K$ of $m$ is either empty, or a point, or a strictly convex body.
\med

{\noi\rm(e)} \ 
The boundary of the core, $\partial K$, is included in $m^*$.
Moreover, $\partial K$ and $m^*$ coincide if and only if the functions
$R_l$ and $R_r$ are nonnegative on $\S^1$.
}
{\noi\sl Remarks.} 
\med\\
1. 
In the context of convex floating bodies, the statement (d) above is already known,
also in arbirtrary dimension,
see {\sl e.g.} Prop. 1 (iv) and (v) in \cit{w} or Theorem 3 in \cit{mr}.
\med\\
2. 
In the case where $R$ vanishes without changing of sign, $m^*$ may,
or may not, be of class $\C^1$.
\med

\pr
{\noi\rm(a)} \ 
Let $s$ denote the curvilinear abscissa on $m$ from some starting point, say $m(0)$. 
We have $s(\q)=\ds\int_0^\q R_l(\tau)d\tau=\ds\int_0^\q R_r(\tau)d\tau$. 
Locally, if $R_l(\q_0)$ and $R_r(\q_0)$ are nonzero, then 
%for $\q$ close to $\q_0$ 
the tangent of $m$ at $m(\q_0)$ is $\vec u(\q_0)$ and we have
$$
\Big(\frac{d\vec u}{ds}\Big)_{l/r}(\q_0)
=\frac{d\vec u}{d\q}\;\Big(\frac{d\q}{ds}\Big)_{l/r}(\q_0)
=\frac1{R_{l/r}(\q_0)}\,\vec u\,'(\q_0).
$$
{\noi\rm(b)} \ 
If $R_l$ and $R_r$ do not vanish, then $s$ is a homeomorphism
(with inverse homeomorphism denoted $\q$ for convenience) and 
$m\big(\q(s)\big)=m(0)+\ds\int_0^s\vec u\big(\q(t)\big)dt$, with
$\vec u$ and $\q$ continuous, {\sl i.e.,} $m\circ\q$ is of class $\C^1$.
\med

{\noi\rm(c)} \ 
We may assume, without loss of generality, that $R_r(\q_1)>0$ and 
$R_l(\q_2)<0$.
If $\q_1=\q_2$ then $m(\q_1)$ is a cusp ({\sl i.e.}, a point with one
half-tangent), hence $m^*$ is not $\C^1$ at $m(\q_1)$.
In the sequel we assume $\q_1<\q_2$.
Consider $E=\{\q\in[\q_1,\q_2]\tq
R_r(\q)>0\}$. Let $\q_3=\sup E$. Since
$0>R_l(\q_2)=\ds\lim_{\q\to\q_2^-}R_l(\q)=\ds\lim_{\q\to\q_2^-}R_r(\q)$,
we have $\q\notin E$ for $\q<\q_2$, $\q$ close enough to $\q_2$.
In the same manner, we have $\q\in E$ if $\q>\q_1$, $\q$ 
close enough to $\q_1$. Therefore, we have $\q_1<\q_3<\q_2$.

Now two cases may occur. 
If there exists $\q_4>\q_3$, such that $R$ vanishes on the whole interval
$\,]\q_3,\q_4[\,$, then $m$ is constant on $[\q_3,\q_4]$. 
We assume $\q_4 \leq \theta_2$ maximal with this property. 
By contradiction, if $m^*$ is $\C^1$ at the point 
$m(\q_3)=m(\q_4)$ then necessarily we have $\q_4=\q_3+\pi$, a contradiction
with $\theta_1<\theta_3\leq\theta_4 \leq \theta_1+\pi$.

In the other case, for all $\delta>0$ there exists $\q\in\,]\q_3,\q_3+\d[\,$
such that $R_r(\q)<0$, and we obtain that $m(\q_3)$ is a cusp.

\med

Before the proofs of (d) and (e), we first establish two lemmas.
%In the rest of the proof,  the core of $m$ is simply denoted by $K$.
%
\lem{l1}{
The interior of the core $K$ of $m$ coincides with the intersection of all
open half-planes $\nt\big(D^+(\q)\big),\;\q\in\S^1$:
\eq j{
\nt(K)=\bigcap_{\q\in\S^1}\nt\big(D^+(\q)\big).
}
Furthermore, we have $\partial K\subset\bigcup_{\q\in\S^1}D(\q)$.
}
\pr
The inclusion $\subseteq$ in~\rf j is evident: For each $\q$ we have
$K\subset D^+(\q)$, hence $\nt(K)\subset\nt\big(D^+(\q)\big)$.
Conversely, given $x\in K$, the map $\q\mapsto\dist\big(x,D(\q)\big)$ 
is continuous.
Since $\S^1$ is compact, if $x$ is in the interior of $D^+(\q)$ for all
$\q\in\S^1$, then this map has a minimum $\rho>0$, 
and the disc of center $x$ and radius $\rho$ is included in $K$.
This proves~\rf j.

If $x\in\partial K$, then $x$ is in every closed half-plane $D^+(\q)$
but not in every open one by~\rf j, hence $x$ has to be on (at least)
one of the lines $D(\q)$.
\ep
\lem{l2}{
If $\q_1\in\S^1$ and $z$ are such that $z\in\partial K\cap D(\q_1)$
then $z=m(\q_1)$.
}
\pr
For small $\eps\neq0$, positive or negative, we have
$$
m(\q_1+\eps)-m(\q_1)=\ds\int_{\q_1}^{\q_1+\eps}R(\q)\vec u(\q)d\q.
$$
By boundedness of $R_l$ and $R_r$, and by continuity of $\vec u$, 
we deduce that there exists $r$ equal to $R_l(\q_1)$ or $R_r(\q_1)$
such that
\eq n{
m(\q_1+\eps)-m(\q_1)=r\vec u(\q_1)\eps+o(\eps).
}
It follows that $D(\q_1)$ and $D(\q_1+\eps)$ cross at a distance $\O(\eps)$
from $m(\q_1)$, see Figure~\reff{f3} below.
If $z$ were different from $m(\q_1)$ then, for $\eps$ small, either negative
or positive depending on the relative positions of $z$ and $m(\q_1)$, 
we would have $z\in\nt\big(D^-(\q_1+\eps)\big)$, hence $z\notin K$,
a contradiction.
\ep

{\noi\sl Now we return to the proof of Proposition~\reff{p1.1}.}
\med

{\noi\rm(d)} \ 
By contradiction, assume that $\partial K$ contains some segment $[x,y]$
with $x\neq y$ and take $z\in\,]x,y[\,$ arbitrarily.
By Lemma~\reff{l1}, there exists $\q_1\in\S^1$ such that $z\in D(\q_1)$.
Since both $x$ and $y$ belong to $K\subset D^+(\q_1)$, the line $D(\q_1)$
contains both $x$ and $y$, {\sl i.e.,} $\q_1$ is the direction
of $\pm\vecc{xy}$.
By Lemma~\reff{l2}, we deduce that $z=m(\q_1)$. 
Hence we proved that {\sl any} $z\in\,]x,y[\,$ coincides with $z=m(\q_1)$, 
where $\q_1$ is one of the directions $\pm\vecc{xy}$, which is impossible.

\med

{\noi\rm(e)} \ 
The first assertion follows directly from Lemmas~\reff{l1} and~\reff{l2}.
For the second one, first we proceed by contradiction and we assume
$R_l(\q_1)<0$ or $R_r(\q_1)<0$ for some $\q_1\in\S^1$,
say $R_r(\q_1)<0$.
Then by~\rf n, for $\eps>0$ small enough, we have 
$m(\q_1+\eps)=m(\q_1)+\eps R_r(\q_1)\vec u(\q_1)+o(\eps)$, hence
$m(\q_1)\in\nt\big(D^-(\q_1+\eps)\big)$, see Figure~\reff{f3}.
It follows that $m(\q_1)\notin K$.
\begin{figure}[htb]
\begin{center}
\begin{picture}(350,50)
\put(45,25){\vector(1,0){225}}
\put(150,25){\circle*{3}}
\put(226,25){\circle*{3}}
\put(59,6){\vector(3,1){140}}
\put(100,19.50){\circle*{3}}
\put(145,10){$m(\q_1)$}
\put(82,3){$m(\q_1+\eps)$}
\put(205,45){$D(\q_1+\eps)$}
\put(272,16){$D(\q_1)$}
\put(225,12){$z$}
\end{picture}
\end{center}
\vspace{-3mm}
\caption{{\small Proof of Lemma~\reff{l2} and Proposition~\reff{p1.1}~(d)}}
\vspace{-5pt}
\lb{f3}
\end{figure}
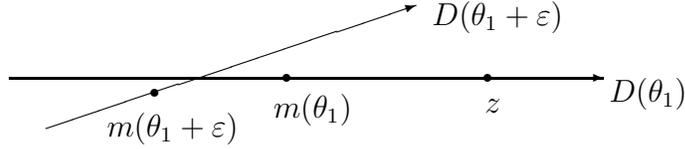

Conversely, if $R_l$ and $R_r$ are nonnegative, take $\q_0\in\S^1$ arbitrarily.
Then, for all $\q\in[\q_0,\q_0+\pi]$, we have
$$
m(\q)=m(\q_0)+\ds\int_{\q_0}^\q\vec m'(\tau)d\tau
=m(\q_0)+\ds\int_{\q_0}^\q R_l(\tau)\vec u(\tau)d\tau.
$$
Therefore we obtain
$$
\langle\vecc{m(\q_0)m(\q)},\vec u\,'(\q_0)\rangle
=\ds\int_{\q_0}^\q R_l(\tau)
\langle\vec u(\tau),\vec u\,'(\q_0)\rangle d\tau\geq0 ,
$$
since $\langle\vec u(\tau),\vec u\,'(\q_0)\rangle =\sin(\q-\q_0)\geq0$. 
In the same manner, we have for all $\q\in[\q_0-\pi,\q_0]$,
$\langle\vecc{m(\q_0)m(\q)},\vec u\,'(\q_0)\rangle\geq0$,
hence the whole curve $m^*$ is in the half-plane $D^+(\q_0)$.
This holds for all $\q_0\in\S^1$, hence $m^*\subset K$.
Finally, for each $\q\in\S^1$, since $m(\q)\in D(\q)$ and 
$K\subset D^+(\q)$, $m(\q)$ cannot be in $\nt(K)$, 
hence $m(\q)\in\partial K$.
\ep
%
%%%%%%%%%%%%%%%%%%%%%%%%%%%%%%%%%%%%%%%%%%%%%%%%%%%%%%%%%%%%%
%
\sec{3.}{Dependence with respect to $\q$}

We begin this section with some known results; we give the proofs
for the sake of completeness of the article. See Section~\reff{2.},
especially Figure~\reff{f1}, for the notation.
We recall that $m$ is the midpoint of the chord $bc$
and that $v_{l|r}=\langle\vec m_{l|r}',\vec u\rangle$.
\propo{p1}{
{\rm(a)} \ 
The curve $m$ is the envelope of the family $\{\De(\a,\q),\q\in\S^1\}$,
{\sl i.e.,} 
\com{http://dictionary.cambridge.org/fr/dictionnaire/britannique/
on-the-one-hand-on-the-other-hand}
%on the one hand 
the vector products 
$\vec m'_l\wedge\vec u$ and $\vec m'_r\wedge\vec u$ 
vanish identically.
%, and on the other hand we have for every open interval
%$I\subset\S^1$ and every continuous %left- and right-differentiable 
%curve $C:I\to K$ such that $C(\q)\in\De(\a,\q)$ for all $\q\in I$, 
%if $\vec C_l'\wedge\vec u=\vec C_r'\wedge\vec u=0$
%identically on $I$ then 
%$C$ is the restriction of $m(\a,\cdot)$ to $I$.
\smallskip

{\noi\rm(b)} \ 
If $b$ and $c$ are regular points of $\partial K$ then
\eq v{
\vec m'=v\vec u~\mbox{ and }~v=\tfrac h2(\cotan\b+\cotan\g).
}
In the case where $b$ or $c$ (or both) is a corner of $\partial K$,
we have $\vec m'_l=v_l\vec u$, $\vec m'_r=v_r\vec u$, 
$v_l=\tfrac h2(\cotan\b_l+\cotan\g_l)$, and 
$v_r=\tfrac h2(\cotan\b_r+\cotan\g_r)$.
\smallskip

{\noi\rm(c)} \ 
If $b$ and $c$ are regular then $v$ is the signed radius of curvature
of the curve $m$.
If $b$ or $c$ (or both) is a corner then $v_l$, resp. $v_r$, is the signed
radius of curvature on the left, resp. on the right, of $m$. 
}
Statement (a) is attributed to M.~M.~Day~\cit{d} by S.~Tabachnikov.
Formula~\rf v appears in a similar form in~\cit{gk}.
\smallskip

\pr
%(a) \
We fix $\a\in\,]0,1[$ and $\q\in\S^1$; we will not always indicate the
dependence in $\a$ and $\q$ of the functions $b,c,m,\vec u$, etc. 
Let $\eps>0$, set $M(\eps)=\De(\a,\q)\cap\De(\a,\q+\eps)$
(see Figure~\reff{f2}), and
consider the curvilinear triangles 
$$
T_b(\eps)=\De(\a,\q)^-\cap\De(\a,\q+\eps)^+\cap K\;,\quad 
T_c(\eps)=\De(\a,\q)^+\cap\De(\a,\q+\eps)^-\cap K.
$$

\begin{figure}[htb]

\begin{center}
\raisebox{0cm}{\epsfysize5.5cm\epsfbox{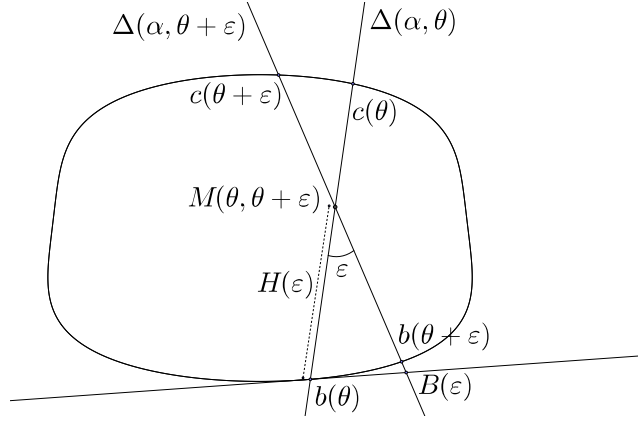}}
%\begin{picture}(300,100)
%\put(83,3){\line(3,1){200}}
%\put(180,50){\oval(150,70)}
%\end{picture}
\end{center}
\vspace{-5mm}
\caption{{\small Proof of Proposition~\reff{p1}}}
\vspace{-5pt}
\lb{f2}
\end{figure}
\bigskip

We have $|\De(\a,\q)^-\cap K|=|\De(\a,\q+\eps)^-\cap K|=\a$, hence
$$0=|T_b(\eps)|-|T_c(\eps)|=
\tfrac12\eps\big(\|b-M(\eps)\|^2-\|c-M(\eps)\|^2\big)
+\o(\eps).
$$
As a consequence, we obtain $\ds\lim_{\eps\to0^+}M(\eps)=m$.
The case $\eps<0$ is similar.

Now we prove
\eq B{
\vec b_r'=h\,\cotan\b_r\,\vec u-h\vec u\,'.
}
For $\eps>0$, let $B(\eps)$ denote the intersection of the right-tangent
to $\partial K$ at $b$ with the line $\De(\a,\q+\eps)$, see Figure~\reff{f2}.
By definition of right-tangent, we have $B(\eps)=b(\q+\eps)+\o(\eps)$.
Set $H(\eps)=\|b-M(\eps)\|$, $H_1(\eps)=\langle B(\eps)-b,\vec u\rangle$,
$H_2(\eps)=\langle M(\eps)-B(\eps),\vec u\rangle$, and 
$H_3(\eps)=\langle b-B(\eps),\vec u\,'\rangle$.
We obtain the following linear system in $H_1,H_2,H_3$
$$
H(\eps)=H_1(\eps)+H_2(\eps),\quad 
H_3(\eps)=H_1(\eps)\tan\b_r=H_2(\eps)\tan\eps,
$$
giving
$$
B(\eps)=b-\frac{\tan\b_r\tan\eps}{\tan\b_r+\tan\eps}H(\eps)\vec u+
\frac{\tan\eps}{\tan\b_r+\tan\eps}H(\eps)\vec u\,'.
$$
As a consequence, we obtain
$B(\eps)=b-\eps H(\eps)\vec u+\eps\cotan\b_rH(\eps)\vec u\,'+\o(\eps)$.
Since $H(\eps)=h+\o(1)$, it follows that
$b(\q+\eps)=B(\eps)+\o(\eps)=
b-\eps h\vec u+\eps\cotan\b_rh\vec u\,'+\o(\eps)$, yielding~\rf B.
Similarly, we have $\vec b_l'=h\,\cotan\b_l\,\vec u-h\vec u\,'$,
$\vec c_l'=h\,\cotan\g_l\,\vec u+h\vec u\,'$, and 
$\vec c_r'=h\,\cotan\g_r\,\vec u+h\vec u\,'$.
Statements (a) and (b) now follow from $m=\frac12(b+c)$.
Since the expression of $v$ given by~\rf v is a ruled function,
statement (c) follows directly from Proposition~\reff{p1.1}~(a).
\ep
\coro{c3}{
For all $\q_0\in\S^1$ and all $\a\in\,]0,1[\,$, we have
\eq u{
\lim_{\q\to\q_0^-}v_l(\q)=\lim_{\q\to\q_0^-}v_r(\q)=v_l(\q_0)
\mbox{ and }
\lim_{\q\to\q_0^+}v_l(\q)=\lim_{\q\to\q_0^+}v_r(\q)=v_r(\q_0).
}
In other words, in the sense of the Pompeiu-Hausdorff distance, we have
\eq V{
\lim_{\q\to\q_0^-}V(\q)=\{v_l(\q_0)\}
\mbox{ and }
\lim_{\q\to\q_0^+}V(\q)=\{v_r(\q_0)\}.
}
Therefore, the function $V$ is  upper semi-continuous (for the inclusion)
with respect to $\q$.
}
\pr
Immediate, using~\rf t, \rf{3b},~\rf v, and the continuity of the cotangent function.
\ep

The next statement will be used both in Sections~\reff{4.} and \reff{5.}.
\coro{c4}{
For $\q_1,\q_2\in\S^1$, with $\q_1\leq\q_2\leq\q_1+\pi$, we have
\eq c{
\conv\big(V(\a,\q_1)\cup V(\a,\q_2)\big)
\subseteq\bigcup_{\q_1\leq\q\leq\q_2}V(\a,\q).
}
}
\pr
Let $v$ be in the above convex hull. If $v\in V(\a,\q_1)\cup V(\a,\q_2)$,
there is nothing to prove.
Otherwise assume, without loss of generality, that 
$\max V(\a,\q_1)<\min V(\a,\q_2)$ and
set $\q_0=\sup\{\q\geq\q_1\tq\max V(\a,\q)<v\}$.
We have $\q_1\leq\q_0<\q_2$; for all $\q\in\,]\q_0,\q_2[$, 
$\max V(\a,\q)\geq v$; and there exists a convergent sequence
$\{\q_n\}_{n\in\N}$ tending to $\q_0$ with $\max V(\a,\q_n)<v$.
By~\rf u, we obtain  $\max V(\a,\q_0)=v$, hence $v\in V(\a,\q_0)$.
\ep
%
%%%%%%%%%%%%%%%%%%%%%%%%%%%%%%%%%%%%%%%%%%%%%%%%%%%%%%%%%%%
%
\sec{4.}{The forwards, the backwards, and the zero sets}
In this section we fix $\a\in\,]0,1[\,$.
Recall that $\q$ is called regular if $\partial K$ is $\C^1$ at $b(\a,\q)$
and $c(\a,\q)$, and singular otherwise.
\df{d2}{
(a) \ 
Let $F(\a)$ be the set of all $\q\in\S^1$ such that either $v(\a,\q)>0$
(if $\q$ is regular),
or $V(\a,\q)\,\cap\,]0,+\infty[\,\neq\emptyset$ (if $\q$ is singular);
we call it the {\sl forwards set}.
\smallskip

\noi(b) \ 
Similarly,  let $B(\a)$ be the set of all $\q\in\S^1$ such that either
$v(\a,\q)<0$ (if $\q$ is regular),
or $V(\a,\q)\,\cap\,]-\infty,0[\,\neq\emptyset$ (if $\q$ is singular);
we call it the {\sl backwards set}.
\smallskip

\noi(c) \ 
Finally, let $Z(\a)$ denote the set of all $\q\in\S^1$ such that either
$v(\a,\q)=0$ (if $\q$ is regular), or $V(\a,\q)$ contains $0$
(if $\q$ is singular); we call it the {\sl zero set}.
}
By the symmetry~\rf W, we have (with the notation
$[\q_1,\q_2]=\{\q\in\S^1\tq\q_1\leq\q\leq\q_2\}$ for $\q_1,\q_2\in\S^1$)
\eq z{
\q\in Z(\a)\ssi\q+\pi\in Z(1-\a)\mbox{ and }
\q\in B(\a)\ssi\q+\pi\in F(1-\a).
}
By Proposition~\reff{p1} (b), we have 
$$
\q\in Z(\a)\ssi\pi\in[\b_l+\g_l,\b_r+\g_r]\;(\mbox{or }
\pi\in [\b_r+\g_r,\b_l+\g_l]).
$$
Since $\b_l\geq\b_r$ and $\g_r\leq\g_l$, by \rf b, we deduce that,
if $\q\in Z(\a)$ then $b=b(\a,\q)$ and $c=c(\a,\q)$ admit parallel 
supporting lines of $K$.
If one of the points $b, c$ or both is regular, the converse is
also true; however, it can occur that
$\b_r+\g_l\leq\pi\leq\b_l+\g_r$ but $\pi$ does not belong to
$[\b_l+\g_l,\b_r+\g_r]$ (or to $[\b_r+\g_r,\b_l+\g_l]$), 
and then $V=[v_l,v_r]$ (or $[v_r,v_l]$) does not contain $0$.

In the case where $\partial K$ is $\C^1$, by Proposition~\reff{p1},
$\q$ belongs to $B(\a)$ if and only if 
$\b(\a,\q)+\g(\a,\q)>\pi$, {\sl i.e.}, there exists a triangle $T=abc$,
with one edge equal to the chord $bc$
(with $b=b(\a,\q),\,c=c(\a,\q)$), the other two edges, $ab$ and $ac$,
not crossing the interior of $K$,
and which contains an amount $1-\a$ of $K$: $|T\cap K|=(1-\a)|K|$.
In the case where $\partial K$ is not $\C^1$, this latter condition is
necessary but not always sufficient to have $\q\in B(\a)$.
However, we will see in Section~\reff{5.} that this condition implies
$\q\in B(\a')$ for all $\a'>\a$.
\med

Let us observe $Z(\a)$ and $B(\a)$, denoted $Z$ and $B$ here.
These sets can be very complicated, even if $\partial K$ is $\C^1$.
In Section~\reff{7.} we present a construction which, to any prescribed
closed subset $C$ of $\S^1$,
associates a $\C^1$ convex body $K$ such that $Z\big(\tfrac12\big)$
coincides with $C$, up to countably many isolated points.
We give next a brief description in the simplest cases.

If $\q_0$ is an isolated point of $Z$  which does not belong to the closure
of $B$, then $m(\q_0)$ is a point of zero
curvature  of $m$, but $m$ is still $\C^1$ at $\q_0$.

If $\q_1<\q_2<\q_1+\pi$ and $[\q_1,\q_2]\subset Z$ is an isolated connected
component of $Z$, then $m$ has a corner; {\sl i.e.,} it is not $\C^1$ but
has two half-tangents: a left-tangent oriented either by $\vec u(\q_1)$ 
if $]\q_1-\d[\,\subset F$ for small $\d>0$, or  by $\vec u(\q_1+\pi)$ if 
$]\q_1-\d[\,\subset B$, and a 
right-tangent oriented either by $\vec u(\q_2)$ or by $\vec u(\q_2+\pi)$.
Moreover, in this situation, $m(\a,\q_1)$ is a local center of symmetry of
$\partial K$: the arc of $\partial K$ in the sector
$\De(\a,\q_1)^+\cap\De(\a,\q_2)^-$ is symmetric to that in
$\De(\a,\q_1)^-\cap\De(\a,\q_2)^+$.

If $\q_0$ is an isolated point of $Z$ and is the endpoint of both a segment
$]\q_1,\q_0[$ of $F$ and a segment $]\q_0,\q_2[$ of $B$,
then $m$ has a cusp at $m(\q_0)$.
\propo{p4}{
Let $\a\in\,]0,1[$.
\smallskip

{\noi\rm(a) \ }
The set $Z(\a)$ is closed  in $\S^1$.
\smallskip

{\noi\rm(b) \ }
If  $\partial K$ is $\C^1$ then the sets $F(\a)$ and $B(\a)$ are open in
$\S^1$ and $F(\a)$, $B(\a)$, and $Z(\a)$ form a partition of \ $\S^1$.
\smallskip

{\noi\rm(c) \ }
In the general case, the three sets $F(\a)\cap B(\a)$, $\partial F(\a)$,
and $\partial B(\a)$ are subsets of $Z(\a)$.
}
{\noi\sl Remark.} 
If $\partial K$ is not $\C^1$ then, in general, there exists $\a\in\,]0,1[\,$
such that $B(\a)$ and $F(\a)$ are not open.
Actually, if $c$ is a corner of $\partial K$, with half-tangents of directions
$\q_1,\q_2$, and $b\in\partial K$ is such that all supporting lines to
$K$ at $b$ have directions in $\,]\q_1+\pi,\q_2+\pi[\,$, then the $\a$-section
(for some $\a$) passing through $b$ and $c$ has a direction 
$\q_0\in F(\a)\cap B(\a)$, 
but such that any $\q<\q_0$, $\q$ close enough to $\q_0$, does not belong to
$B(\a)$ and any $\q>\q_0$, $\q$ close enough to $\q_0$, does not belong 
to $F(\a)$.
\med

\pr
(b) \ 
If $\partial K$ is $\C^1$ then the functions $\b$ and $\g$ are continuous
with respect to $\q$, hence also is $v$,
so $Z=v^{-1}(\{0\})$ is a closed subset of $\S^1$,
and $F=v^{-1}(\,]-\infty,0[\,)$ and  
$B=v^{-1}(\,]0,+\infty[\,)$ are open. Furthermore they do not intersect
and their union is $\S^1$.
\med

\noi(a) \ Let $\{\q_n\}_{n\in\N}$ be a sequence converging to  $\q_0$, with
$\q_n\in Z(\a)$, {\sl i.e.,} $0\in V(\a,\q_n)$.
Then a subsequence $\{\q_{n_k}\}_{k\in\N}$ exists such that  $\q_{n_k}$
tends either to $\q_0^-$ or  to $\q_0^+$, say  to $\q_0^-$. 
By~\rf V we have
$\ds\lim_{k\to+\infty}V(\a,\q_{n_k})=\big\{v_l(\a,\q_0)\big\}$,
hence $v_l(\a,\q_0)=0$, and $\q_0\in Z(\a)$.
\med

\noi(c) \ Let $\q_0\in\partial F(\a)$; then there exist sequences
$\{\q_n\}_{n\in\N}$ in $F(\a)$ and $\{\q'_n\}_{n\in\N}$ in
$\S^1\setminus F(\a)$, both converging to $\q_0$. 
This means that $V(\q_n)\cap\,]0,+\infty[\,\neq\emptyset$ and
$V(\q'_n)\subset\,]-\infty,0]$ for every $n\in\N$. 
By~\rf c, there exists $\q''_n\in[\q_n,\q'_n]$ (or $\q''_n\in[\q'_n,\q_n]$)
such that $0\in V(\q''_n)$ for each $n$, {\sl i.e.,} $\q''_n\in Z(\a)$.
Since the sequence $\{\q''_n\}_{n\in\N}$ tends to $\q_0$ and $Z(\a)$
is closed, we obtain  $\q_0\in Z(\a)$.
The proof for $\partial B(\a)\subseteq Z(\a)$ is similar.
The proof for $F(\a)\cap B(\a)\subseteq Z(\a)$ is obvious.
\ep
%
%%%%%%%%%%%%%%%%%%%%%%%%%%%%%%%%%%%%%%%%%%%%%%%%%%%%%%%%%%%%%%%%%%
%
\sec{5.}{Dependence with respect to $\a$}
The functions $b$ and $c$ are left- and right-differentiable with respect
to $\a$; one finds, {\sl e.g.},
$\veccc{\lp\tfrac{\partial b}{\partial\a}\rp_l}
=\tfrac1{2h}(\vec u\,'+\cotan\b_r\,\vec u)$.
However, we will not use their differentiability in $\a$, but only their
monotonicity. 
Since the function $\cotan$ is decreasing on $\,]0,\pi[\,$, by~\rf g we
immediately obtain the following statement, whose proof is omitted.
\propo{p1bis}{
The functions $v,v_l$ and $v_r$ are nonincreasing in $\a$.
More precisely, we have 
\eq m{
\max V(\a',\q)\leq\min V(\a,\q)\mbox{ for all }\q\in\S^1
\mbox{ and all }0<\a<\a'<1.
}
}
Recall the sets $F(\a)$, $B(\a)$, and $Z(\a)$ introduced in 
Definition~\reff{d2}. Let 
\begin{eqnarray*}
I_F&=&\{\a\in\,]0,1[\,\tq F(\a)\neq\emptyset\},\\
I_B&=&\{\a\in\,]0,1[\,\tq B(\a)\neq\emptyset\},\\
I_Z&=&\{\a\in\,]0,1[\,\tq Z(\a)\neq\emptyset\}.
\end{eqnarray*}
The symmetry~\rf z implies
\eq i{
\a\in I_F\ssi1-\a\in I_B\mbox{ and }\a\in I_Z\ssi1-\a\in I_Z.
}
Set $\a_B=\inf I_B$ and $\a_Z=\inf I_Z$. (We do not consider $\inf I_F$,
which is  always equal to $0$, see below).
\theo{t2}{
{\noi\rm(a)} \ 
We have $I_B=\,]\a_B,1[$ (and hence $I_F=\,]0,1-\a_B[$ by~\rf i).
If $\partial K$ is $\C^1$, then  $I_Z=[\a_Z,1-\a_Z]$; otherwise $I_Z$
is one of the intervals $[\a_Z,1-\a_Z]$ or $]\a_Z,1-\a_Z[$.
\med

{\noi\rm(b)} \ We have $\a_Z\leq\a_B\leq\frac12$. 
Moreover, if $\a_Z\neq\a_B$ then $\partial K$ contains two parallel segments.
\med

{\noi\rm(c)} \ 
We have the following equivalences:
\smallskip

{\rm i.} \ $\a_B=\tfrac12$ if and only if $K$ is symmetric.
\smallskip

{\rm ii.} \ $\a_Z=0$ if and only if $\partial K$ contains a segment whose
endpoints admit two parallel supporting lines to $K$.
\smallskip

{\rm iii.} \ $\a_B=0$ if and only if $\partial K$ contains a segment whose
endpoints admit two parallel supporting lines to $K$, one of which
intersecting $\partial K$ at only one point.
}
{\noi\sl Remarks.} 
1. Notice the change of behaviour of $Z$ and $B$:
In Proposition~\reff{p4}, $Z(\a)$ is always closed, whereas $B(\a)$ may be
not open if $\partial K$ is not $\C^1$; here $I_B$ is always open, whereas
$I_Z$ may be not closed if $\partial K$ is not $\C^1$.
\med

\noi2. An example of a planar convex body such that $I_Z$ is open is the 
quadrilateral $OICJ$ in Figure~\reff{f7}, Section 8.9, for which one
finds $I_Z=\,\big]\frac1{2c},1-\frac1{2c}\big[\,$.
\med

\pr
\noi(a) \
By Proposition~\reff{p1bis},
if $0<\a<\a'<1$ then $B(\a)\subseteq B(\a')$. 
It follows that $I_B=\,]\a_B,1[\,$  or $I_B=[\a_B,1[$. 
We now prove that $\a_B\notin I_B$.

If $\a\in I_B$, then there exists $\q\in\S^1$ such that, say, $v_l(\a,\q)<0$
(the case $v_r(\a,\q)<0$ is similar).
Let $\a'<\a$ and $\q'<\q$ be such that $b(\a',\q')=b(\a,\q)$; then by~\rf u
\footnote{
In fact \rf u is stated for $\a$ fixed, but its proof uses only~\rf t
and~\rf v --- which can easily be adapted to
our situation --- and the continuity of the function $\cotan$.
}%
, both $v_l(\a',\q')$ and $v_r(\a',\q')$ tend to $v_l(\a,\q)$ as $\q'\to\q$
and $\a'\to\a$, hence are negative for $(\a',\q')$ close enough to
$(\a,\q)$. This shows that $\a'\in I_B$. As a consequence, $I_B=\,]\a_B,1[\,$.

We now prove that $I_Z$ is convex. Let $\a_1<\a_2\in I_Z$; {\sl i.e.},
there exist $\q_1$ and $\q_2$ in $\S^1$ such that 
$0\in V(\a_1,\q_1)\cap V(\a_2,\q_2)$.
We may assume, without loss of generality, that $\q_1\leq\q_2\leq\q_1+\pi$.
Let $\a\in\,]\a_1,\a_2[$.
By~\rf m, we have $V(\a,\q_1)\cap\R^-\neq\emptyset$ and
$V(\a,\q_2)\cap\R^+\neq\emptyset$ hence, by~\rf c,
we have $0\in V(\a,\q)$ for some $\q\in[\q_1,\q_2]$
(in the case $\q_2=\q_1+\pi$, both intervals $[\q_1,\q_2]$ and
$[\q_2,\q_1+2\pi]$ suit). This shows that $\a\in I_Z$.
\med

\noi(b) \
If $\a$ belongs neither to $I_F$ nor to $I_B$, then we necessarily have
$V(\a,\q)=\{0\}$ for all $\q\in\S^1$, hence the function $m(\a,\cdot)$
is constant. This implies $\a=\frac12$ and $K$ symmetric. As a consequence,
we have $\a_B\leq\frac12$.

We prove $\frac12\in I_Z$, yielding $I_Z\neq\emptyset$.
Let $\q_0\in\S^1$ arbitrary. By~\rf W,
$0\in\conv\big(V\big(\frac12,\q_0\big)\cup V\big(\frac12,\q_0+\pi\big)\big)$,
hence by~\rf c $0\in V\big(\frac12,\q\big)$ for some $\q\in[\q_0,\q_0+\pi]$,
hence $Z\big(\tfrac12\big)\neq\emptyset$.

We now prove $\a_Z\leq\a_B$.
If $\a_B=\frac12$, we are done; otherwise, let $\a\in I_B$, $\a\leq\frac12$.
Then $v_l(\a,\q)$ or $v_r(\a,\q)$ is negative for some $\q\in\S^1$, 
say $v_l(\a,\q) <0$.
From~\rf w and~\rf m it follows that $v_l(\a,\q+\pi)=-v_l(1-\a,\q)\geq-v_l(\a,\q)>0$,
hence by continuity there exists $\q_1\in\,]\q,\q+\pi[$ such that
$v_l(\a,\q_1)=0$, so $\q_1\in Z(\a)$, and $\a\in\ I_Z$.
\med

If moreover $\a_Z<\a_B$, then consider $\a_Z<\a<\a'<\a_B\leq\frac12$,
hence $\a,\a'\in I_Z\setminus I_B$. 
Let $\q\in Z(\a)$, $b=b(\a,\q)$, $c=c(\a,\q)$, and let $D_b$, $D_c$ be two
parallel supporting lines to $K$ at $b$ and $c$ respectively,
see Figure~\reff{f4}. 
\begin{figure}[htb]
\begin{center}
\raisebox{0cm}{\epsfysize5.5cm\epsfbox{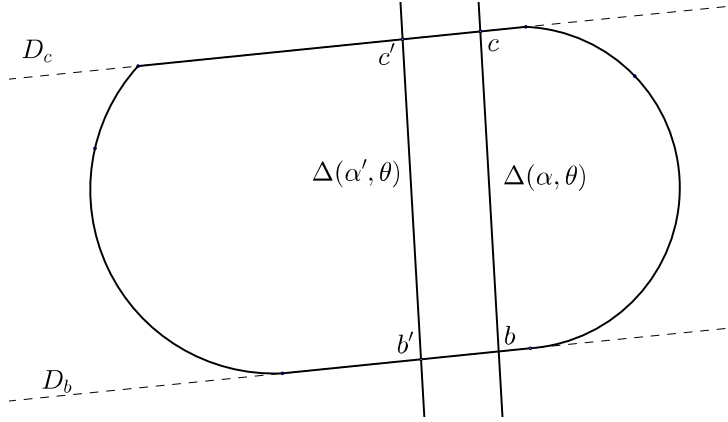}}
%\begin{picture}(300,100)
%\put(83,3){\line(3,1){200}}
%\put(180,50){\oval(150,70)}
%\end{picture}
\end{center}
\vspace{-5mm}
\caption{{\small Proof of Theorem~\reff{t2}~(b)}
}
\vspace{-5pt}
\lb{f4}
\end{figure}
Consider now $b'=b(\a',\q)$ (for the same $\q$) and $c'=c(\a',\q)$.
Since $\a<\a'$, $\De(\a',\q)$ is in the interior of
$\De^+(\a,\q)$. Since $B(\a')=\emptyset$, there do not exist supporting lines
to $K$ at $b'$ and $c'$ crossing in $\De^+(\a',\q)$.
As a consequence, $b'$ must lie on $D_b$ and $c'$ on $D_c$, and the segments
$[b',b]$ and $[c',c]$ are on $\partial K$.
\med

\noi(c) i. \ If $\a_B=\frac12$ then by (a) $\frac12\notin I_F\cup I_B$ hence,
as already said in the proof of (b), $K$ is symmetric.
The converse is obvious.
\med

\noi(c) ii. and iii. \ If $\a_Z=0$ then there exists a sequence
$\{\a_n\}_{n\in\N}$ tending to $0$, with $\a_n\in I_Z$;
{\sl i.e.}, such that for any $n$ there exists $\q_n\in\S^1$ with
$0\in V(\a_n,\q_n)$.
By compactness, we may assume without loss of generality that the sequence
$\{\q_n\}_{n\in\N}$ converges to some $\q_0\in\S^1$.
Let $b_n=b(\a_n,\q_n)$ and $c_n=c(\a_n,\q_n)$.
By continuity, the sequences  $\{b_n\}_{n\in\N}$ and  $\{c_n\}_{n\in\N}$
converge to some $b$, resp. $c\in\partial K$.
Since the width of $K$ satisfies $0<w(K)\leq\|b_n-c_n\|$, we have $b\neq c$.
Since $0\in V(\a_n,\q_n)$, there are parallel supporting lines to $K$ at
$b_n$ and $c_n$ for each $n$.
Without loss of generality, we also assume that these lines converge,
yielding two parallel supporting lines at $b$ and $c$,
denoted $D_b$ and $D_c$.
Since the sequence of lines $\big\{\De(\a_n,\q_n)\big\}_{n\in\N}$ tends
to the oriented line through $b$ and $c$, denoted by $(bc)$, and
$|\De^-(\a_n,\q_n)\cap K|=\a_n|K|\to0$, we obtain $|(bc)^-\cap K|=0$,
hence the segment $[b,c]$ is on $\partial K$.

Conversely, if $[b,c]\subset\partial K$ and $b\neq c$ admit two parallel
supporting lines to $K$, denoted by $D_b$ and $D_c$,
let $\q$ denote the direction of $(bc)$;
{\sl i.e.,} $\q$ is such that $\vec{bc}=k\vec u(\q)$ for some $k>0$.
The angles $\b_n=\b\big(\frac1n,\q\big)$ and $\g_n=\g\big(\frac1n,\q\big)$
satisfy $\b_n+\g_n\geq\pi$,
hence $\q\in B\big(\frac1n\big)\cup Z\big(\frac1n\big)$,
so $B\big(\frac1n\big)\cup Z\big(\frac1n\big)\neq\emptyset$,
$\frac1n\in I_B\cup I_Z$, and
$0\leq\min(\a_Z,\a_B)=0$, hence $\a_Z=0$ since $\a_Z\leq\a_B$.

If one of the supporting lines above, say $D_b$, intersects $\partial K$
only at $b$, then $b_n\notin D_b$, yielding $\b_n+\g_n>\pi$,
hence $\frac1n\in I_B$, showing $\a_B=0$.
Conversely, if $\a_B=0$, then the former points $b_n$ and $c_n$ admit
supporting lines which cross in $\De^+(\a_n,\q_n)$ (see comment after
Definition~\reff{d2}).
By contradiction, if both $D_b\cap\partial K$ and $D_c\cap\partial K$
contained more than $b$, resp. $c$,
then for $n$ large enough we would have $b_n\in D_b$ and $c_n\in D_c$.
Then for $m,n$ such that $b_m\in\,]b_n,b[\,$ and $c_m\in\,]c_n,c[\,$,
the only supporting lines at $b_m, c_m$ could be $D_b$ and
$D_c$, which do not cross, a contradiction.
\ep
%
%%%%%%%%%%%%%%%%%%%%%%%%%%%%%%%%%%%%%%%%%%%%%%%%%%%%%%%%%%%%%
%
\sec{6.}{The $\a$-core}
In this section we compare the boundary $\partial K_\a$ with the image of
$m(\a,\cdot)$, denoted by $m^*_\a$, for $\a\in\,]0,1[\,$.
The function of $\q$, $m(\a,\cdot)$, is a tour in the sense of
Definition~\reff{d1}, with $R_l=v_l(\a,\cdot)$ and $R_r=v_r(\a,\cdot)$,
and core $K_\a$. In the case $\a=\frac12$,  $m\big(\frac12,\cdot\big)$
is also a double half-tour.

Then, by Proposition~\reff{p1.1} (b) and (c) and Theorem~\reff{t2} (a),
$m^*_\a$ is of class $\C^1$
if $0<\a<\a_Z$ or  $1-\a_Z<\a<1$,
and is not $\C^1$ if $\a_B<\a<1-\a_B$. In the case $\a=\a_Z$, 
$m^*_\a$ may or may not be $\C^1$. 
However, if $\a_Z\neq\a_B$, we will see that $m^*_\a$ is not $\C^1$ also
for $\a_Z<\a\leq\a_B$.

By Proposition~\reff{p1.1} (d), $K_\a$ is strictly convex
(or one point or empty), 
and by Proposition~\reff{p1.1} (e) we have 
$m^*_\a=\partial K_\a$ if and only if $B(\a)$ is empty.
By Theorem~\reff{t2}~(a), we then have
\eq8{
m^*_\a=\partial K_\a\ssi\a\leq\a_B.
} 
Besides, the function $\a\mapsto K_\a$ is continuous (for the 
Pompeiu-Hausdorff distance on compact sets in the plane)
and decreasing (with respect to inclusion):  if $\a<\a'$, then
$K_{\a'}\subsetneq K_\a$.
Since $K_\a=\emptyset$ if $\a>\frac12$, there exists a value
$\a_K\leq\frac12$ such that  
\begin{itemize}
\item
$K_\a$ is strictly convex with a nonempty interior if $0<\a<\a_K$,  
\item
$K_{\a_K}$ is a single point denoted by $T$, and 
\item
$K_\a$ is empty if $\a_K<\a<1$.
\end{itemize}
By~\rf8, this value $\a_K$ is at least $\a_B$.
It is noticeable that, when $K_{\a_K}$ is a single point, this point
is {\sl not} necessarily the mass center of $K$, see Section~8.9.

We end this section with the following statement, which gathers the
last results.
\propo{t3}{
{\noi\rm(a)} \ If $\a_Z<\a<1-\a_Z$, then $m^*_\a$ is not $\C^1$.
\med

{\noi\rm(b)} \ We have $\frac49\leq\a_K\leq\frac12$, with first equality
if and only if $K$ is a triangle,
and second equality if and only if $K$ is symmetric.
\med

{\noi\rm(c)} \ If $K$ is non-symmetric, then $\a_B<\a_K$
(whereas for $K$ symmetric we have  $\a_B=\a_K=\frac12$).
}
\pr
(a)
It remains to prove that, if $\a_Z<\a\leq\a_B$, then $m_\a^*$ is not $\C^1$.
Assume $\a_Z<\a_B$; by Theorem~\reff{t2}~(b) and its proof,
there exist two parallel segments on $\partial K$, denoted by
$[a,b]$ and $[d,c]$ with $abcd$ in convex position, such that
the line oriented by $\vec{bc}$, denoted by $D_1$, is an $\a_Z$-section of $K$.
Let $D_2$ denote the line oriented by $\vec{ad}$. 
%We have $|D_1^-\cap K|\leq|D_2^+\cap K|$, otherwise the line
%oriented by $\vec{da}$, {\sl i.e.} opposite to $D_2$, would be 
%an $\a$-section with some 
%thus $D_1$ is an $\a_Z$-section of $K$ and 
Then $D_2$ is a $\b$-section for some $\b\geq\a_B$.
For any $\a\in\,]\a_Z,\a_B]$, $\a > \a_Z$ but close to it,
there is an $\a$-section, denoted by 
$\De(\a,\q_1)$, passing through $c$ and crossing
$[a,b]$ and an $\a$-section, denoted by 
$\De(\a,\q_2)$, passing through $b$ 
and crossing $[d,c]$, see Figure~\reff{f5}.
These $\a$-sections intersect at some point 
$P\in\big[\frac12(a+d),\frac12(b+c)\big[\,$.
Then we have $m(\a,\q)=P$ for all $\q\in[\q_1,\q_2]$, 
$\vec{m'_l}(\a,\q_1)$ collinear to $\vec{Pc}=\|\vec{Pc}\|\,\vec u_{\q_1}$,
and $\vec{m'_r}(\a,\q_2)$ collinear to $\vec{bP}=\|\vec{bP}\|\,\vec u_{\q_2}$.
Since $\q_1\neq\q_2\neq\q_1+\pi$, $m(\a,\cdot)$ is not $\C^1$ at $P$.
If $\a$ is not close to $\a_Z$, $\De(\a,\q_1)$ (passing through $c$) will not cross $[a,b]$,
and $\De(\a,\q_2)$ (passing through $b$) will not cross $[d,c]$.
Nevertheless, the argument given above holds as well,
except that $\De(\a,\q_1)$ and $\De(\a,\q_2)$ will not cross on 
$\big[\frac12(a+d),\frac12(b+c\big[$.
\begin{figure}[htb]
\begin{center}
\raisebox{0cm}{\epsfysize5.5cm\epsfbox{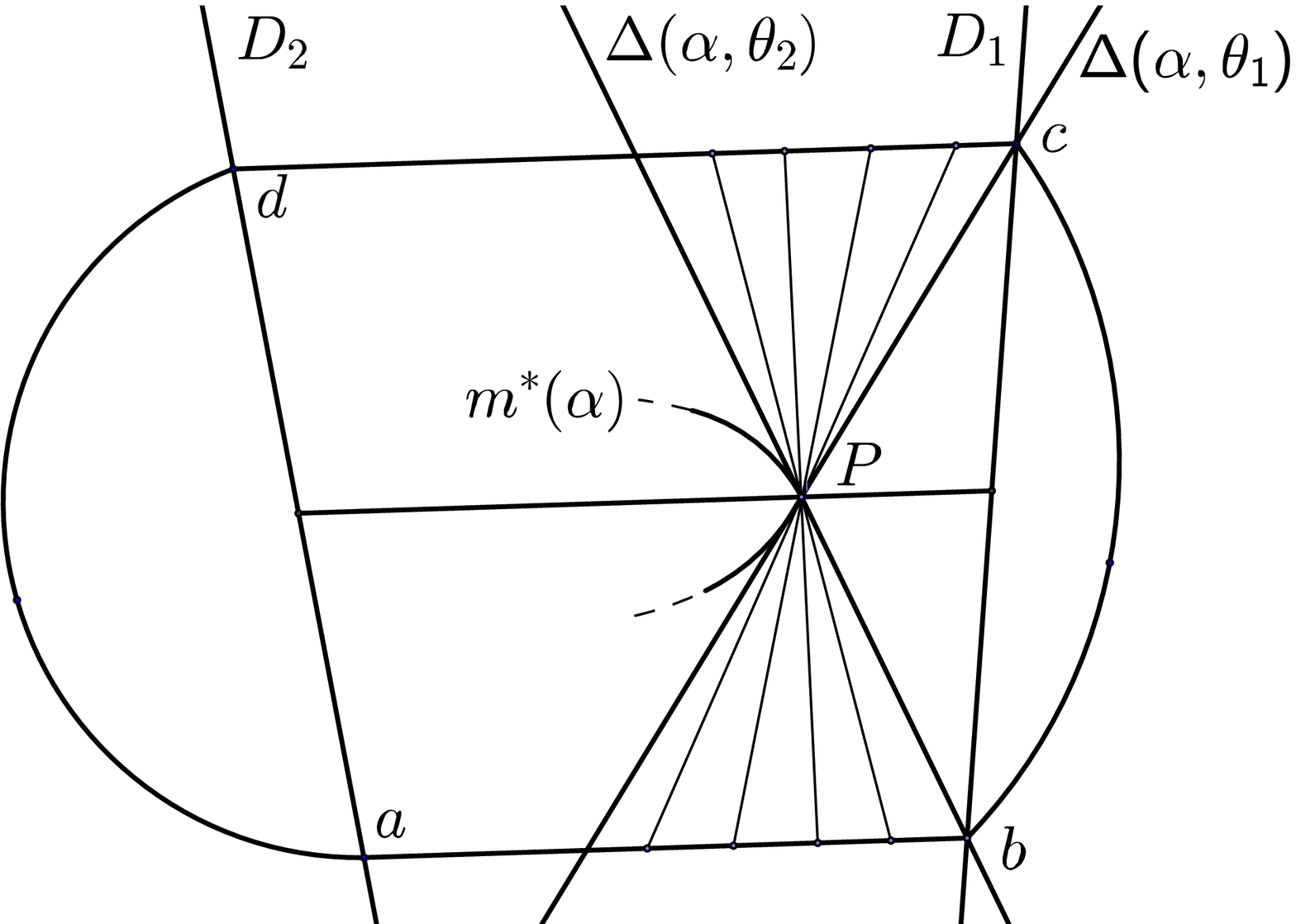}}
%\begin{picture}(300,100)
%\put(83,3){\line(3,1){200}}
%\put(180,50){\oval(150,70)}
%\end{picture}
\end{center}
\vspace{-5mm}
\caption{{\small Proof of Proposition~\reff{t3} (a)}
}
\vspace{-5pt}
\lb{f5}
\end{figure}
\med

\noi(b) 
Suppose $\a_{K}=\frac{1}{2}$. Let $T$ be the unique point of $K_{1/2}$.
With the notation for the proof of Proposition~\reff{p1}~(a),
see Figure~\reff{f2}, since every half-section contains $T$,
we have $M(\q,\q')=T$ for all $\q\neq\q'$.
Now this proof implies that
$$
\left\vert \q-\q'\right\vert (\left\Vert b(\q
)-T\right\Vert -\left\Vert c(\q)-T\right\Vert )=\O\big(\left\vert \q
-\q'\right\vert ^2\big),
$$
for all $\q<\q^{\prime}$. Letting $\q^{\prime}$ tend to $\q$,
this shows that $b(\q)$ and $c(\q)$ are symmetric about $T$.

The fact that $\a_K\geq\frac49$ and that $\a_K=\frac49$ 
if and only if $K$ is a triangle is well known;
see, {\sl e.g.},~\cit n.
\com{j'ai retrouve la ref !!}
It reduces to the following statement which we prove below.
\smallskip

{\sl The mass center $G$ of $K$ is in $\De^+\big(\frac49,\q\big)$
for all $\q$. Moreover, there exists $\q\in\S^1$ such that 
$\De(\frac49,\q)$ contains $G$ if and only if $K$ is a triangle.}
\smallskip

Fix $\q$, let $\De=\De\big(\tfrac49,\q\big)$, and consider the frame 
$\big(b=b(\q),\vec u=\vec u(\q),\vec v=\vec{u'}(\q)\big)$.
We have to prove that $y$-coordinate $G_y$ of $G$ is nonnegative and that
this coordinate vanishes if and only if $K$ is a triangle.
Let $D_b$ and $D_c$ be two supporting lines of $K$ at the points $b$ and
$c=c(\q)$, see Figure~\ref{f9}. 
\begin{figure}[htb]
\begin{center}
\raisebox{0cm}{\epsfysize6.5cm\epsfbox{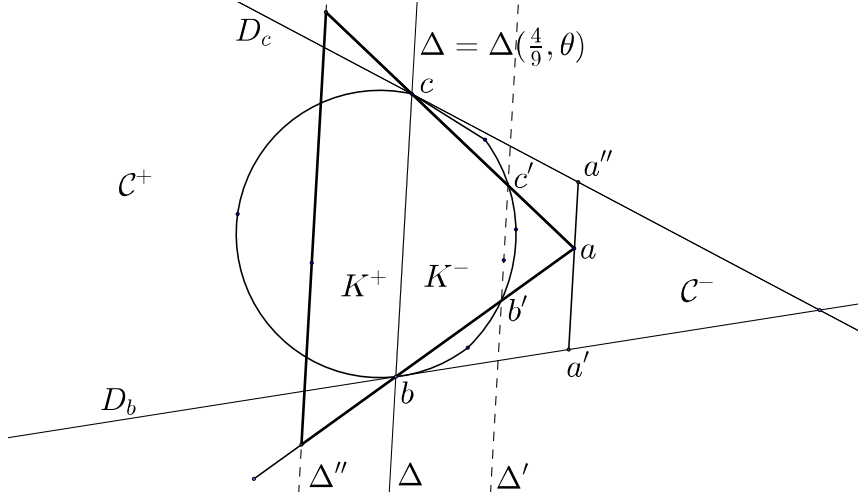}}
%\begin{picture}(300,100)
%\put(83,3){\line(3,1){200}}
%\put(180,50){\oval(150,70)}
%\end{picture}
\end{center}
\vspace{-5mm}
\caption{{\small Proof of Proposition~\reff{t3} (b). The 
triangle $V=L\cup U$ is in bold.}
}
\vspace{-5pt}
\lb{f9}
\end{figure}
These lines together with the segment $[b,c]$ define two convex
sets $\mathcal C^+=\conv\big((\De^+\cap
D_b)\cup(\De^+\cap D_c)\big)$ and $\mathcal C^-=\conv\big((\De^-\cap D_b)
\cup(\De^-\cap D_c)\big)$, one on each side of $\De$.
The set $K^+=K\cap\De^+$ is included in $\mathcal C^+$ and
the set $K^-=K\cap\De^-$ is included in $\mathcal C^-$. 
Let $a'\in D_b\cap\De^-$ and $a''\in D_b\cap\De^-$ be such that
the triangles $\conv(a',b,c)$ and $\conv(a'',b,c)$ have an area equal to
$\frac49|K|$. Then the line $(a'a'')$ is parallel to $\De$ and for every 
$a\in[a',a'']$ the triangle $\conv(a,b,c)$ also has an area equal to
$\frac49|K|$.
For any such $a$, let $b'\in\partial K\cap[a,b[$ and 
$c'\in\partial K\cap[a,c[$. By continuity, one can choose $a$ such that 
the oriented line $\De'=(b'c')$ is parallel to $\De$. For this $a$, let 
$L=\conv(a,b,c)$. We have $K^-\setminus L\subset\De'^+$ and 
$L\setminus K^-\subset\De'^-$, therefore the $y$-coordinate of the 
mass center $G^-$ of $K^-$ is larger than or equal to the 
$y$-coordinate of the mass center of $L$, the equality holding only when 
$K^-$ is a triangle.
Furthermore, since $a\in\mathcal C^-$, the convex set $\mathcal C^+$ is
included in the convex set
$\mathcal C^{\prime+}=\conv\big((\De^+\cap(ab)\cup(\De ^+\cap(ac)\big)$.

Let $U$ be the trapezoid of area $\frac59|K|$ defined by the lines
$(ab)$, $(ac)$, $\De$ and a line $\De''\subset\De^+$ parallel to $\De$.
Since $K^+\setminus U\subset\De''^+$ and $U\setminus K^+\subset\De''^-$,
the $y$-coordinate of the mass center $G^+$ of $K^+$ is larger than or
equal to the $y$-coordinate of the mass center of $U$, the equality 
holding only when $K^+=U$. 
It follows that $G_y$ is larger than or equal to the $y$-coordinate of the
mass center of the union $V$ of $L$ and $U$, with equality only if $K=V$.
Since the mass center of $V$ is on $\De$ we are done.
\med

\noi(c)
Suppose that $\a_B=\a_K$. By~\rf8, we obtain that $m(\a_K,\cdot)$
is constant and equal to $T$. It follows that, for all $\q\in\S^1$,
$T$ is the middle of a chord of direction $\q$, hence $K$ is symmetric 
about $T$.
\ep
%
%\com{doit-on encore ajouter un bilan ? 
%
%$0\leq\a_Z\leq\a_B\leq\a_K\leq\frac12$, ce qu'il se passe entre chaque valeur,
%en separant $K$ sym ou non, un peu comme dans l'intro ?
%\\
%ou bien ecrire un th resumant le tout dans l'intro directement ?
%}
%
%
%
%%%%%%%%%%%%%%%%%%%%%%%%%%%%%%%%%%%%%%%%%%%%%%%%%%%%%%%%
%
\sec{7.}{Miscellaneous results, remarks, and questions}
{\noi\bf8.1 \ }
It is easy to check that, if an $\a$-section crosses two non-parallel
segments of $\partial K$, then the corresponding middle of chord $m(\a)$
lies on an arc of hyperbola asymptotic to the lines extending these segments. 
In particular, if $K$ is a convex polygon, then for all $\a\in\,]0,1[$ the
curve $m(\a,\cdot)$ is entirely made of arcs of hyperbolae.
Of course, not only segments yield arcs of hyperbolae.
One can check for instance that two arcs of hyperbolae on $\partial K$ also
give an arc of hyperbola for $m(\a,\cdot)$, if $\a$ is small enough.
\bigskip

{\noi\bf8.2 \ }
We saw that an envelope of $\a$-sections is a tour in the sense of 
Definition~\reff{d1}.
Conversely, which tours are envelopes of some $\a$-sections $m(\a,K)$,
for which convex bodies $K$ and for which $0<\a<1$?
For example, the astroid on the left of Figure~\reff{f8} given by $\q\mapsto(-\cos^3\q,\sin^3\q)$, {\sl i.e.,} 
with $R(\q)=\frac32\sin(2\q)$, cannot be such an envelop. 
More generally, following a remark by D.~Fuchs and S.~Tabachnikov
in~\cit{ft}, a tour with a common tangent at two different points
cannot be such an envelop, since each point would have to be the
middle of the chord.
\bigskip

{\noi\bf8.3 \ }
The KAM theory applied to dual billiards shows that, if $m$ is a tour
of class $\C^5$ (hence strictly convex), then there exist 
convex bodies $K$ with $\partial K$ arbitrarily close to $m$ such that
$m=m(\a,K)$ (hence for some $\a$ arbitrarily close to $0$) and there exist 
convex bodies $K$ with $\partial K$ arbitrarily close to infinity such that
$m=m(\a,K)$ (hence for some $\a$ arbitrarily close to $\frac12$).
The curves $\partial K$ are invariant torii of  the dual billiard. 
According to E.~Gutkin et A.~Katok \cit{gk}, the works of J.~~Moser~\cit{m}
and R.~Douady~\cit{d} prove that these curves are convex.
Notice the following apparent paradox. By Theorem~\reff{t2}(c)i., if $K$ is
non-symmetric then, for $\a$ close to $\frac12$, $m(\a)$ is non-convex.
This apparently contradicts the above results which imply that, given a
strictly convex $\C^5$ non-symmetric curve, 
there exist values of $\a$ arbitrarily close to $\frac12$ and convex bodies
$K$, necessarily non-symmetric, such that $m=m(\a,K)$. 
However, these curves have envelops $m(\b,K)$ with cusps for some other
$\b\in\,]\a,\frac12[$, although
$\a$ is arbitrarily close to $\frac12$. There is no real contradiction.
\med

Several questions remain open: 
Among tours which present cusps, which ones are envelopes of sections?
For instance,  does it exist a {\sl symmetric} curve with cusps which is the
$\a$-envelope of some (necessarily non-symmetric) convex body?
Also here, does it exist a non-symmetric convex body $K$ having a symmetric
convex envelope $m(K,\a,\cdot)$ for some $\a\in\,]0,1[\,$?
\bigskip

{\noi\bf8.4 \ }
The link between envelopes of $\a$-sections and dual billiards yields
a simple proof for the following (also simple) fact:

The only convex bodies which have a circular envelope (or an elliptic one,
since this is the same question modulo an affine transformation)
are the discs having the same centers as the envelop.
Actually, if the billiard table is a circle, then the billiard map is
integrable
\footnote{A classical conjecture states that only ellipses have a dual
billiard map which is integrable.}.
This means that each orbit remains on a circle centered at the origin,
and is either periodic or dense in this circle, 
depending whether the angle between the two tangents from the starting
point at the table is in $\pi\Q$ or not. 
If the convex body $K$ were not a disc, then its boundary would cross
at least one circle with dense orbits, 
hence would contain at least one dense orbit on the latter circle,
hence the whole circle, a contradiction.
\bigskip

{\noi\bf8.5 \ }
Let $C$ be a closed subset of $\S^1$. We construct below a $\C^1$ convex curve
which is the boundary of a convex body for which $Z\big(\frac12\big)=C\cup S$,
for some countable set $S$.

As convex curve, we start with the unit circle and will deform it.
Since $C$ is closed, $\S^1\setminus C$ is a countable union of open intervals.
Let $]\q_1,\q_2[$ be one of them.
Between $\q_1$ and $\q_2$, we deform the circle into a convex $\C^1$ curve
arbitrarily close to the union of two segments, one $[e^{i\q_2},p]$
tangent to the circle at $e^{i\q_2}$, the other $[e^{i\q_1},p]$, such that
the area remains unchanged.
(For convenience, here we use the notation of complex numbers.)
We do the same symmetrically with respect to the line passing through the
origin and of direction $\frac12(\q_1+\q_2+\pi)$; {\sl i.e.},
we choose a $\C^1$ curve close to $[e^{i(\q_1+\pi)},p']$ tangent to the circle
at $e^{i(\q_1+\pi)}$, and to $[e^{i(\q_2+\pi)},p']$,
with $p'=\overline{p}\,e^{i(\q_1+\q_2+\pi)}$, see Figure~\reff{f6}.

Observe that the curve is no longer (centrally-)symmetric. 
If the curve is chosen with large curvature near the points $p$ and $p'$,
then the interval $]\q_1,\q_2[$ contains only two values of $\q$
for which the half-section cuts the curve in a chord with two parallel
tangents; these chords have one endpoint near $p$, resp. near $p'$, and 
do not contain the origin anymore.

Doing this in any connected component of $\S^1\setminus C$, we obtain a
convex $\C^1$ curve such that any $\q\in C$ has a halving chord with parallel
tangents and all but two values of $\q$ in each connected component of
$\S^1\setminus C$ have a halving chord without parallel tangents.
\begin{figure}[htb]
\begin{center}
%\begin{picture}(300,100)
%\put(83,3){\line(3,1){200}}
%\put(180,50){\oval(150,70)}
%\end{picture}
\raisebox{0cm}{\epsfysize6cm\epsfbox{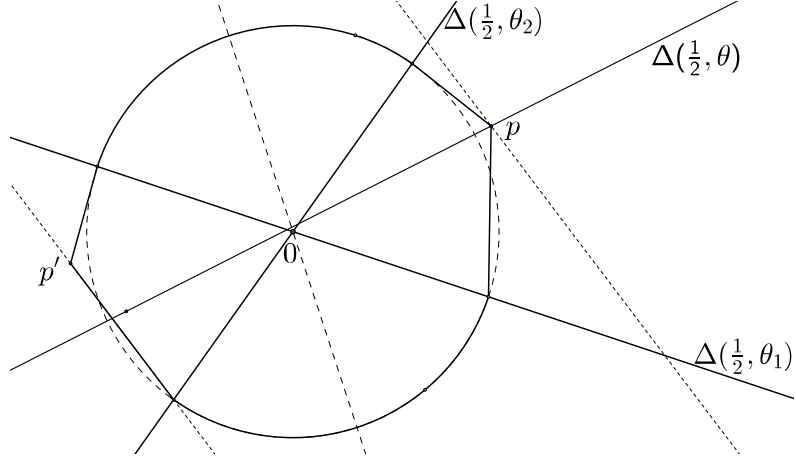}}
\end{center}
\vspace{-5mm}
\caption{{\small A construction of $\C^1$ convex body with (almost)
prescribed zero set}
}
\vspace{-5pt}
\lb{f6}
\end{figure}
\bigskip

{\noi\bf8.6 \ }
It is easy to construct $K\varsubsetneq L$ such that $L_\a\subset K_\a$
for some value of $\a$;
for example, $L$ is an equilateral triangle, $K$ the inscribed circle of $L$,
and $\a\in\big[\frac49,\frac12\big]$.

Our calculation shows that this holds for all
$\a\in\big[\a_1,\frac12\big]$, where $\a_1\approx0.40716$ satisfies
$1-\frac32(1-\sqrt{1-2\a_1})=\cos t$ and $\a_1=\frac1\pi(t-\cos t\sin t)$.
Actually, if the radius of $K$ is $R$ (hence the height of $L$ is $3R$),
$K_\a$ is a disc of radius $R_\a=R\cos t$ such that
$\a=\frac1\pi(t-\cos t\sin t)$.
Then $\partial L_\a$ is made of three arcs of hyperbolae.
In a (non-orthonormal) frame where the triangle has vertices $(0,0)$, 
$(0,1)$, and $(1,0)$, two of these arcs of hyperbolae are 
$xy=\frac\a4$ and $x(1-x-y)=\frac\a4$,
hence they cross at a point of abscissa $x_c=\frac12(1-\sqrt{1-2\a})$)
This proves that the three arcs of $\partial L_\a$ cross at
a distance $R-\frac32R(1-\sqrt{1-2\a})$ from the center.

Is it true that, for every value of $\a\in\,]0,\tfrac12[$ there exists
$K\subsetneq L$ such that $L_\a\subset K_\a$?
(It cannot exist $K,L$, $K\subsetneq L$, independent of $\a$ such that
$L_\a\subset K_\a$ for all small $\a>0$).

On the contrary, does it exist $\a_1>0$ such that for all pairs $K,L$,
$K\subsetneq L$, and all $\a<\a_1$, $L_\a\not\subset K_\a$?
Is $\a_1\approx0.40716$ optimal with this respect? {\sl I.e.,} is it
optimal when $K$ is a disc and $L$ a triangle?
Or, instead, is it optimal when $K$ is an affinely regular hexagon
inscribed in $L$, a triangle?
\bigskip

{\noi\bf8.7 \ }
The following questions have been asked by Jin-ichi Itoh, whom we would
like to thank for his interest in our work.
Let $w(K)$ denote the width of a convex body $K$, $\oslash(K)$ its diameter,
$r(K)$ its inradius, and $R(K)$ its circumradius.
\med

Do we have $\frac{w(K)}{\oslash(K)}\leq\frac{w(K_\a)}{\oslash(K_\a)}$
for all $\a<\frac12$?
\med

Do we have $\frac{r(K)}{R(K)}\leq\frac{r(K_\a)}{R(K_\a)}$ for all
$\a<\frac12$?
\med

Both answers are ``no'', as shown by the following counter-examples.

For the first question, if $K$ is the unit Reuleaux triangle,
$\oslash(K)=w(K)=1$, and $\a>0$ is small enough, then one can check that
$\oslash(K_\a)=1-\O\big(\a^{2/3}\big)$, whereas
$w(K_\a)<1-\a^{1/2}<\oslash(K_\a)$.
\med

For the second question, consider two small arcs of the circle of center
0 and radius 2 near the $x$-axis, two small arcs of the circle of center
0 and radius 1 near the $y$-axis, and choose for $K$ the convex hull of
the union of these four arcs.

In this manner we have $\frac{r(K)}{R(K)}=\frac12$. However, we have
$r(K_\a)=\cos t_1$, where $t_1$ is such that 
$
\a=\tfrac1{|K|}(t_1-\sin t_1\cos t_1)=\tfrac{2}{3|K|}t_1^3+\O(t_1^5),
$
and $R(K_\a)=2\cos t_2$, where $t_2$ is such that 
$
\a=\tfrac4{|K|}(t_2-\sin t_2\cos t_2)=
\tfrac{8}{3|K|}t_2^3+\O(t_2^5),
$
hence $t_2<t_1$, so $\frac{r(K_\a)}{R(K_\a)}<\frac12$.

Whether there is a constant $k>0$ such that
$\frac{w(K)}{\oslash(K)}\leq k\frac{w(K_\a)}{\oslash(K_\a)}$
for all convex bodies $K$ and all $\a$, idem for $\frac{r(K)}{R(K)}$,
and which constant is optimal, seems to be another interesting question.
\bigskip

{\noi\bf8.8 \ }
There is a tight link between the region where appear cusps,
as described in~\cit{ft}, and the regions $M_2$ and $M_3$
described by T. Zamfirescu in~\cit z. 
Actually, let $K$ be a planar convex body and, for each $\q\in\S^1$,
consider the so-called {\sl midcurve}, {\it i.e.}, the locus of midpoints
of all chords of direction $\q$.
If $\partial K$ is $\C^1$ and strictly convex then the family of
all these midcurves for all $\q\in\S^1$ is a continuous family in
the sense of Gr\"unbaum~\cit g, for which general results
of~\cit z are at our disposal.

In our situation, the regions $M_2$ and $M_3$ of~\cit z are the loci
of points of $\nt(K)$ which are middles of at least $2$, resp. $3$,
different chords.
If $K$ is symmetric then $M_2=M_3=\{g(K)\}$, the mass center of $K$.
If $K$ is not symmetric then it seems that we have 
\eq6{
\nt(M_2)=M_3=\nt\Big(\!\!\bigcup_{\a\in\,]0,\frac12[}(m_\a^*\setminus\partial K_\a)\Big).
}
Already~\rf6 can be verified when $K$ is a quadrilateral. For instance,
the polygonal curve of all cusps described in Figure 11.15 of~\cit{ft} 
is indeed the boundary of $M_3$.

Besides, when $K$ is a quadrilateral, a detailed analysis shows that
the unique point $T$ of $K_{\a_K}$ is in $M_3$
and that the three chords bisected by $T$ are all $\a_K$-sections.
This property of $T$ is likely to be true in the general case.
Conversely, is $T$ the {\sl only} point of
$M_3$ bisecting three $\a$-sections with the same $\a$?
This has been checked in the case of quadrilaterals.
\bigskip

{\noi\bf8.9 \ }
We present here an explicit example showing that the unique point $T$ of
$K_{\a_K}$ is not necessarily the mass center $G$ of $K$.
The strategy of proof is to use a quadrilateral and to show that the
three chords bisected by $G$ are $\a$-sections for different values of $\a$,
showing that $T\neq G$ by Subsection~8.8.

Given $c>1$, let $K=OIJC$ denote the quadrilateral $\conv(O,I,J,C)$,
with $O=(0,0),\,I=(1,0),\,J=(0,1),\,C=(c,c)$.
The area of $K$ is $|K|=\frac12\,\|\veccc{IJ}\|\,\|\veccc{OC}\|=c$ and its
mass center is $G=\big(\frac{2c+1}6,\frac{2c+1}6\big)$
(the midpoint of the segment determined by the centers of mass of the
triangles $OJC$ and $OIC$).
The point $G$ is the midpoint of exactly three chords, $AB$, $EF$,
and $E'F'$, with $A,\,B$ on the line of equation $x+y=\frac{2c+1}3$,
$E\in OJ$, $F\in IC$, and $E'\in OI$, $F'\in JC$, see Figure~\reff{f7}.
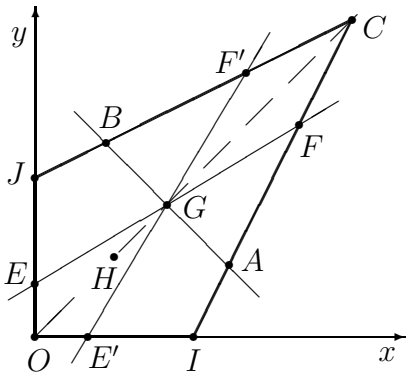
\begin{figure}[htb]
\begin{center}
\begin{picture}(140,120)
\put(0,-.3){\line(1,0){60}}\put(0,0){\vector(1,0){140}}
\put(0,.3){\line(1,0){60}}
\put(-.3,0){\line(0,1){60}}\put(0,0){\line(0,1){60}}
\put(.3,0){\vector(0,1){125}}
\multiput(32,32)(16,16){6}{{\line(1,1){10}}}\put(0,0){\line(1,1){15}}
\put(0,59.7){\line(2,1){120}}\put(0,60){\line(2,1){120}}
\put(0,60.3){\line(2,1){120}}
\put(59.7,0){\line(1,2){60}}\put(60,0){\line(1,2){60}}
\put(60.3,0){\line(1,2){60}}
\put(-10,14){\line(5,3){125}}
\put(14,-10){\line(3,5){75}}
\put(15,85){\line(1,-1){70}}
\put(130,-9){$x$}\put(-9,112){$y$}
\put(0,0){\circle*{3}}\put(-3,-13){$O$}
\put(60,0){\circle*{3}}\put(57,-13){$I$}
\put(0,60){\circle*{3}}\put(-11,58){$J$}
\put(120,120){\circle*{3}}\put(124,113){$C$}
%\put(45,45){\circle*{3}}\put(36,45){$T$}
\put(30,30){\circle*{3}}\put(56,45){$G$}
\put(50,50){\circle*{3}}\put(21,18){$H$}
\put(0,20){\circle*{3}}\put(-12,20){$E$}
\put(100,80){\circle*{3}}\put(100,68){$F$}
\put(20,0){\circle*{3}}\put(20,-12){$E'$}
\put(80,100){\circle*{3}}\put(68,100){$F'$}
\put(27,73.5){\circle*{3}}\put(24,79){$B$}
\put(73.5,27){\circle*{3}}\put(78,25){$A$}
\end{picture}
\end{center}
%\vspace{-3mm}
\caption{{\small The $\a$-center is not the mass center; here $c=2$}}
\vspace{-5pt}
\lb{f7}
\end{figure}
These chords divide $K$ in proportion $(\a_{AB},1-\a_{AB})$, resp.
$(\a_{EF},1-\a_{EF})$ and $(\a_{E'F'},1-\a_{E'F'})$
(by convention, $0<\a_{XY}\leq\frac12$).
By symmetry, we have $\a_{EF}=\a_{E'F'}$, but there is no reason that
$\a_{EF}=\a_{AB}$; as we will see, this is indeed not the case.

The triangle $ABC$ is the image of $IJC$ by the dilation of center $C$
sending the point $H=\big(\frac12,\frac12\big)$ to $G$,
hence of factor $\frac{\|\veccc{GC}\|}{\|\veccc{HC}\|}=\frac{4c-1}{3(2c-1)}$.
With $|IJC|=c-\frac12$, this gives $|ABC|=\frac{(4c-1)^2}{18(2c-1)}$, so
$\a_{AB}=\frac{(4c-1)^2}{18c(2c-1)}$.

Since $F$ is on the line $IC$ of equation $c(x-1)=(c-1)y$,
one finds $F=\big(\frac{2c+1}3,\frac{2c}3\big)$, hence $E=\big(0,\frac13\big)$.
The area of the triangle $EIF$ is 
$\frac12\det(\veccc{EI}\;\veccc{EF})=\frac12\left|
\begin{matrix} 1&\frac{2c+1}3\\
-\frac13&\frac{2c-1}3
\end{matrix}\right|=\frac{4c-1}9$, 
hence the area of the quadrilateral $OIFE$ is
$|OIFE|=|OIE|+|EIF|=\frac16+\frac{4c-1}9=\frac{8c+1}{18}$.
This gives $\a_{EF}=\frac{8c+1}{18c}\neq\a_{AB}$, hence $G\neq T$.
\com{on pourrait aussi calculer T explicitement}
\med

More generally, given a planar convex body $G$, we consider all chords 
the midpoint of which coincides to the mass center $G$ of $K$. 
It is known~\cit v that there are at least three different such chords; 
these chords are $\a$-sections for different values of $\a$ 
ranging from some $\a_{\min{}}(K)$ to some $\a_{\max{}}(K)$.
Then the quotient $\frac{\a_{\min{}}(K)}{\a_{\max{}}(K)}$ measures in 
some sense the asymmetry of the body $K$. It is equal to $1$ if $K$ is 
(centrally) symmetric, or if $K$ has another symmetry which ensures that $G$ 
is the only affine-invariant point, {\sl e.g.} for $K$ a regular polygon, 
but it differs from $1$ in general.
One can see that the minimum of this quotient is achieved for at least one 
affine class of convex bodies, and it would be worthy to determine the shape
of these bodies, and to compute the corresponding minimum.
In the case of all quadrilaterals, using Maple, we found that the minimum is 
attained precisely for a quadrilateral of our previous family,
for $c=\frac74$ and with $\frac{\a_{\min{}}(K)}{\a_{\max{}}(K)}=\frac{24}{25}$.
The same questions can be asked in higher dimensions, replacing the
midpoint of a chord by the mass center of an $\a$-section.
\bigskip

{\noi\bf8.10 \ }
We end this section with our main conjecture.
\conj{c1}{
For any  convex bodies $K, L$ with $K\subset L$, and any
$\a\in\,]0,\frac12[\,$, there exists an $\a$-section of $L$ which is a
$\b$-section of $K$ for some $\b\leq\a$.
}
This conjecture has been recently proven in~\cit{fm} in the case of planar 
convex bodies.
Another natural related question, the answer of which turns out to be negative,
is the following.
\med

For $K\subset L$ convex bodies, does there always exist a half-section
$\De$ of $L$ such that
$$
|K|\lng(\De\cap L)\leq|L|\lng(\De\cap K)\;?
$$
For a counter-example, consider for $K$ a thin pentagon of width varying
from $\eps$ on each side to $2\eps$ in the middle, placed at the basis 
of $L$, a triangle of heigth $1$, see Figure~\reff{f18}.
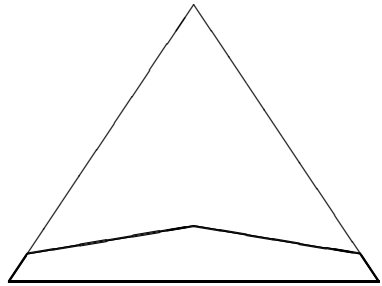
\begin{figure}[htb]
\begin{center}
\begin{picture}(140,110)
\setlength{\unitlength}{0.7pt}
\put(0,-.3){\line(1,0){200}}\put(0,0){\line(1,0){200}}
\put(0,.3){\line(1,0){200}}
\put(0,0){\line(2,3){100}}
\put(100,150){\line(2,-3){100}}
\put(0,-.3){\line(2,3){10}}\put(0,0){\line(2,3){10}}
\put(0,.3){\line(2,3){10}}
\put(10,14.7){\line(6,1){90}}\put(10,15){\line(6,1){90}}
\put(10,15.3){\line(6,1){90}}
\put(100,29.7){\line(6,-1){90}}
\put(100,30){\line(6,-1){90}}\put(100,30.3){\line(6,-1){90}}
\put(190,14.7){\line(2,-3){10}}\put(190,15){\line(2,-3){10}}
\put(190,15.3){\line(2,-3){10}}
%\put(100,16.367){\line(6,-1){100}}
%\put(100,16.667){\line(6,-1){100}}\put(100,16.967){\line(6,-1){100}}
\end{picture}
\end{center}
\vspace{-3mm}
\caption{{\small
In bold, the body $K$; in thin, the body $L$; here $\eps=\frac1{10}$}
}
\vspace{-5pt}
\lb{f18}
\end{figure}
A detailed analysis shows that, for every half-section $\De$ of $L$,
one has $|K|\lng(\De\cap L)>|L|\lng(\De\cap K)$.

\

{\noi\bf Acknowledgements.} The authors thank Theodor Hangan and Tudor
Zamfirescu for fruitful discussions and for having pointed to them several
references.
They are also indebted to the anonymous referee who drew to their attention
the reference~\cit{mr}. 
The third author thanks the {\sl Universit� de Haute Alsace}
for a one-month grant and its hospitality, 
and acknowledges partial support from the Romanian National Authority for
Scientific Research, CNCS-UEFISCDI, grant PN-II-ID-PCE-2011-3-0533.
%\footnotesize

\vspace{1cm}

{\small 
\noi Addresses of the authors:
\med

\noi Nicolas Chevallier and Augustin Fruchard\\
Laboratoire de Math�matiques, Informatique et Applications\\
Facult� des Sciences et Techniques\\
Universit� de Haute Alsace\\
2 rue des Fr�res Lumi�re\\
68093 Mulhouse cedex, FRANCE\med\\
E-mails: {\tt Nicolas.Chevallier@uha.fr,\ 
Augustin.Fruchard@uha.fr}

\

\noi Costin V\^{\i}lcu\\
Simion Stoilow Institute of Mathematics of the Romanian Academy,\\
P.O. Box 1-764, Bucharest 70700, ROMANIA\med\\
E-mail: {\tt Costin.Vilcu@imar.ro}

\end{document}